\documentclass{amsart}
\usepackage{graphicx}
\usepackage{amssymb,amsmath}

\theoremstyle{plain}  
\newtheorem{thm}{Theorem}

\newtheorem{lemma}[thm]{Lemma}  

\newtheorem{prop}[thm]{Proposition}

\newtheorem{rem}[thm]{Remark}

\numberwithin{equation}{section}

\newcommand{\al}{\alpha}
\newcommand{\C}{\mathbb C}

\newcommand{\Q}{\mathbb Q}
\newcommand{\R}{\mathbb R}
\newcommand{\Z}{\mathbb Z}

\newcommand{\PSLC}{PSL(2,\mathbb C)}
\newcommand{\SLC}{SL(2,\mathbb C)}

\newcommand{\ve}{\varepsilon}

\begin{document}
\title{Finite surgeries on three-tangle pretzel knots}
\author[D. Futer, M. Ishikawa, Y. Kabaya, T. Mattman, and K. Shimokawa]{David Futer, Masaharu Ishikawa, Yuichi Kabaya, Thomas W.\ Mattman, and Koya Shimokawa}
\address{Mathematics Department, Temple University, Philadelphia, PA 19122, USA}
\email{dfuter@math.temple.edu}
\address{Mathematical Institute, Tohoku University,
Sendai, 980-8578, Japan}
\email{ishikawa@math.tohoku.ac.jp}
\address{Department of Mathematics, Tokyo Institute of Technology, 2-12-1 Oh-okayama, Meguro-ku, Tokyo 152-8551, Japan}
\email{kabaya@math.titech.ac.jp}
\address{Department of Mathematics and Statistics,  California State University, Chico,  Chico CA 95929-0525, USA}
\email{TMattman@CSUChico.edu}
\address{Department of Mathematics,  Saitama University,  Saitama 338-8570, Japan}
\email{kshimoka@rimath.saitama-u.ac.jp}

\thanks{The second author is supported by MEXT, Grant-in-Aid for Young Scientists (B) (No. 19740029). 
This research occured during a visit of the fourth author to Saitama University. He would like to thank the Math Department for their warm hospitality. 
The fifth author is  partially supported by the Ministry of Education, Culture, Sports, Science and Technology, Grant-in-Aid for Scientific Research (C) 18540069}
\begin{abstract}
We classify Dehn surgeries on $(p,q,r)$ pretzel knots that result in a
manifold of finite fundamental group. The only hyperbolic pretzel knots
that admit non-trivial finite surgeries are $(-2,3,7)$ and $(-2,3,9)$.
Agol and Lackenby's 6-theorem reduces the argument to knots with small
indices $p,q,r$. We treat these using the Culler-Shalen norm of 
the $\SLC$-character variety. In particular, we introduce
new techniques for demonstrating that boundary slopes are 
detected by the character variety.
\end{abstract}

\dedicatory{Dedicated to Professor Akio Kawauchi on the occasion of his
60th birthday.}

\maketitle

In \cite{m2} Mattman showed that if a hyperbolic $(p,q,r)$ pretzel knot $K$ admits
a non-trivial finite Dehn surgery of slope $s$ (i.e., a Dehn surgery that results
in a manifold of finite fundamental group) then either
\begin{itemize}
\item $K = (-2,3,7)$ and $s = 17$, $18$, or $19$, 
\item $K = (-2,3,9)$ and $s = 22$ or $23$, or
\item $K = (-2,p,q)$ where $p$ and $q$ are odd and 
$5 \leq p \leq q$.
\end{itemize}
In the current paper we complete the classification by proving
\begin{thm}
\label{thmmain}
Let $K$ be a $(-2,p,q)$ pretzel knot with $p$, $q$ odd and
$5 \leq p \leq q$. Then $K$ admits no non-trivial finite surgery.
\end{thm}

Using the work of Agol~\cite{a} and Lackenby~\cite{l},
candidates for finite surgery correspond to curves of length at most six in the maximal cusp of $S^3 \setminus K$.
%
If $7 \leq p \leq q$, we will argue that only five slopes for the $(-2,p,q)$ pretzel knot have length
six or less: the meridian and the four integral surgeries $2(p+q)-1$, $2(p+q)$, $2(p+q) + 1$, and $2(p+q)+2$.
If $p = 5$ and $q \geq 11$, a similar argument leaves seven candidates,
the meridian and the six integral slopes between $2(5+q) -2$ and 
$2(5+q)+3$. 

We will treat the remaining knots, $(-2,5,5)$, $(-2,5,7)$, and $(-2,5,9)$,  using the Culler-Shalen norm (for example, see
\cite{bz2,cgls}). For a hyperbolic knot in $S^3$, this is a norm
$\| \cdot \|$ on the vector space $H_1(\partial M; \R)$. We can 
identify a Dehn surgery slope $s \in \Q \cup \{ \frac10 \}$ with a class 
$\gamma_s \in H_1(\partial M ; \Z)$. If $s$ is a finite slope that 
is not a boundary slope, 
the finite surgery theorem~\cite{bz2} shows that $s$ is integral
or half-integral and $\| \gamma_s \| \leq \max \{ 2S, S + 8 \}$
where $S = \min \{ \| \gamma \| : 0 \neq \gamma \in H_1(\partial M; \Z) \}$
is the minimal norm. This makes the Culler-Shalen norm an effective
tool for the study of finite surgery slopes.

The Culler-Shalen norm is intimately related to the set of boundary
slopes. An essential surface in the knot complement $M$ will meet $\partial M$ in a (possibly empty) set of parallel curves. The slope
represented by this set of curves is known as a {\it boundary slope}. 
For a pretzel knot, these slopes are determined
by the algorithm of Hatcher and Oertel~\cite{ho}. Given the list of
boundary classes $\{ \beta_j : 1 \leq j \leq N \}$, the norm is determined by an associated set of non-negative integers $a_j$:
$$
\| \gamma \| = 2 \sum_{j = 1}^{N} a_j \Delta( \gamma, \beta_j).
$$
Here $\Delta(.,.)$ denotes the minimal geometric intersection number.

In particular, if a boundary class $\beta_j$ is detected by the character variety, then the corresponding coefficient $a_j$ is positive.  To describe the notion of ``detection,'' recall that,  following Culler and Shalen~\cite{cs}, we can use $\SLC$-representations of the knot group to construct essential surfaces in the complement $M$. The construction uses an ideal point of a curve in the $\SLC$-character variety to induce a non-trivial action of the knot group on a Bass-Serre tree. This action in turn yields an embedded essential surface in $M$. We will say that the boundary slope associated to a surface so constructed is {\it detected} by the character variety. As it will be enough for us to argue that certain $a_j$ are positive, an important part of our proof is the introduction of new techniques to show that a boundary slope is detected.

We remark that Ichihara and Jong have recently announced an independent proof of Theorem~\ref{thmmain} using Heegaard Floer homology \cite{ii}. 
Indeed, they go further and classify finite and cyclic surgeries for
all Montesinos knots.
Also, with an argument based on Khovanov Homology, Watson~\cite{wa} has just shown that the $(-2,p,p)$ pretzel knot (for $p$ odd, $5 \leq p \leq 25$) admits no non-trivial finite surgeries.

Our paper is organised as follows. As above, our proof of 
Theorem~\ref{thmmain} breaks into two cases. In Section~\ref{sec6thm}, we
use the 6-theorem to handle the case where
$p \geq 7$ as well as the case where $p = 5$ and $q \geq 11$ and prove that none of these knots admit non-trivial finite surgeries.
This leaves the three knots $(-2,5,5)$, $(-2,5,7)$, and $(-2,5,9)$,
which we treat using the Culler-Shalen norm. 
In Section~\ref{sec255}, we show that  
the $(-2,5,5)$ pretzel knot admits no non-trivial finite 
surgeries by introducing new techniques for detecting boundary
slopes that generalise the method introduced by Kabaya~\cite{k}
(based on earlier work of Yoshida).
In Section~\ref{sec2579}, we give some observations concerning
detection of the boundary slopes $2(p+q)$ for the $(-2,5,5)$ and $(-2,5,7)$ pretzel knots using techniques pioneered by Ohtsuki~\cite{o1,o2}.
In consequence, we conclude that the $(-2,5,7)$ pretzel knot admits no non-trivial finite surgeries. We conclude Section~\ref{sec2579} by proving
the same assertion for the $(-2,5,9)$ pretzel knot. Thus, in each section, we show that the knots under consideration admit no non-trivial finite slopes. Taken together, this proves Theorem~\ref{thmmain}.

In the next section we collect
some general results that will be used throughout the paper.

\section{\label{seclem}%
Lemmas}

In this section let $K$ denote a $(-2,p,q)$ pretzel knot
with $p,q$ odd and $5 \leq p \leq q$. We collect several facts about
finite slopes of these knots. We begin with arguments that apply to
all slopes. We next look at arguments specific to even integer slopes
and those that apply to the slopes $2(p+q) \pm 1$.

For a pair of slopes $\frac{a}{b},\frac{c}{d} \in \Q \cup \{\frac10\}$, the distance is defined to be
$\Delta( \frac{a}{b}, \frac{c}{d} ) = |ad-bc|$. This is equivalent to the minimal geometric intersection number of curves representing these two slopes. Agol~\cite{a} and Lackenby~\cite{l} independently showed that any pair of exceptional slopes on a one-cusped hyperbolic manifold lie within distance $10$ of each other. Very recently, Lackenby and Meyerhoff \cite{lm} improved the bound from $10$ to $8$; we will not need this improvement.

Since $2(p+q)$ is an exceptional, toroidal surgery slope
of the $(-2,p,q)$ pretzel knot $K$~\cite{w}, it follows from Agol and Lackenby's work \cite{a, l} that
any other exceptional slope $s$ is within distance $10$ of $2(p+q)$.
\begin{lemma}
\label{lem610}
Let $p,q$ be odd and $5 \leq p \leq q$. If $s$ is a finite slope
of the $(-2,p,q)$ pretzel knot, then
$\Delta( s, 2(p+q)) \leq 10$. 
\end{lemma}

The next set of lemmas relate to the Culler-Shalen norm for a hyperbolic knot in $S^3$; \cite{s} is a good reference. Note that $K$ is hyperbolic~\cite{kaw}.
We will be working with the total norm, which is given by summing the Culler--Shalen norm over all components
of the character variety. The following lemma, which is proved in Section~5.1 of~\cite{m1} (c.f., \cite[Proposition~1.3]{m2}), gives
the minimal total norm.
\begin{lemma}
\label{lemS}%
Let $p,q$ be odd with $5 \leq p \leq q$.
The minimal total norm of the $(-2,p,q)$ pretzel knot is
$S = 2pq - 3(p+q)$.
\end{lemma}

The next lemma summarises several consequences of the finite surgery theorem~\cite{bz2}.
\begin{lemma}
\label{lembz}
Let $p,q$ be odd with $5 \leq p \leq q$ and let 
 $s = \frac{a}{b}$ be a finite surgery slope of the 
$(-2,p,q)$ pretzel knot. Then,
either $s = \frac10$ or else $s$ is integral (i.e., $b=1$) or half-integral
($b = 2$). Moreover, if $\gamma_s$ is the class of a curve representing 
slope $s$ and $s$ is not a boundary slope, then $\| \gamma_s \| \leq S+8$ unless $s$ is an even integer,
in which case $\| \gamma_s \| \leq 2S$.
\end{lemma}
Note that $M(\frac10) = S^3$ is known as {\it trivial surgery} and 
$\frac10$ is always a finite slope. Our goal in this paper is
to show that it is the only finite slope of $K$.

In earlier work, Ishikawa, Mattman, and Shimokawa showed that there is a relationship between
the positions of finite slopes and boundary slopes~\cite{ims}. 
Recall that a surface $F \subset M$ is called a \emph{semi-fibre} if its complement $M \setminus F$ is a (possibly twisted) $I$--bundle. (Under this definition, a fibre in a fibration of $M$ is one example of a semi-fibre.) 
A boundary slope is \emph{strict} if it is not the slope of a semi-fibre of $S^3 \setminus K$.
For a knot in $S^3$, if $r$ is a boundary slope that is not strict, then $r = 0$.

\begin{lemma}
\label{lemims}
Let $K$ be a hyperbolic knot in $S^3$ and $\| \cdot \|$ the associated
total Culler-Shalen norm. Let $S$ denote the minimal total norm and suppose that $S \geq 8$. If $s=\frac{a}{b}$ is a finite slope of $K$, then there is a strict boundary slope $r$ with $|s-r| \leq 2/b$ and the inequality is
strict unless $K$ has exactly two strict boundary slopes detected by the character variety.
\end{lemma}

\begin{proof}
In the proof of Corollary 3(2) of \cite{ims}, it was shown
that $|s-r| \leq 3/b$ using $S \geq 4$. If we instead assume
$S \geq 8$, the same argument shows that $|s-r| \leq 2/b$, as required.
\end{proof}


Finally, we observe the following immediate consequence of \cite[Theorem 2.0.3]{cgls}.
Recall that a knot is {\em small} if there is no closed essential surface in its complement.
In particular, the $(-2,p,q)$ pretzel knot $K$ is small~\cite{ot}.

\begin{lemma}
\label{lemfinbdy}
Let $\beta$ be a boundary slope for a small knot $K$ in $S^3$. Then $\beta$ is not a finite slope.
\end{lemma}

\subsection{Even surgeries}

In this subsection we will see that even integer surgeries of the
$(-2,p,q)$ pretzel knot are generally not finite.

Using the Wirtinger presentation~\cite{r}, the fundamental group
of the $(-2,p,q)$ pretzel knot is (c.f., \cite{t})
\begin{eqnarray*}
\pi_1(M) = \langle x,y,z &\mid&
(zx)^{(p-1)/2}z(zx)^{(1-p)/2} = (yx)^{-(q+1)/2}y(yx)^{(q+1)/2},\\
&&(yz^{-1})^{-1}y(yz^{-1}) = (yx)^{(1-q)/2}x(yx)^{(q-1)/2},\\
&&(yz^{-1})^{-1}z(yz^{-1}) = (zx)^{(p+1)/2}x(zx)^{-(p+1)/2}
\rangle .
\end{eqnarray*}
There is redundancy in the relations as any one is a consequence 
of the other two.

For an integral surgery slope $s$, the group of $M(s)$ is given
by adding the relator $x^sl$ where 
$$l = x^{-2(p+q)}(yx)^{(q-1)/2}(yz^{-1})^{-1}(yx)^{(q+1)/2}(zx)^{(p-1)/2}
(yz^{-1})(zx)^{(p+1)/2},$$
represents the longitude, i.e., $\pi_1(M(s)) = \pi_1(M) / \langle x^s l \rangle.$

\begin{lemma}
\label{lemeven}%
Let $p,q$ be odd with $5 \leq p \leq q$. Assume further that if $p = 5$,
then $q \geq 11$. Let $s$ be an even integer and let $M$ be the complement
of the $(-2,p,q)$ pretzel knot. Then $\pi_1(M(s))$ is 
not finite.
\end{lemma}

\begin{proof}
We will argue that  $\pi_1(M(s))$ projects onto the group
that Coxeter~\cite{c} calls $(2,p,q;2)$. Edjvet~\cite{e} has shown that this group is infinite under our hypotheses on $p$ and $q$.

Assume that $s$ is even.
Adding the relators $x^2, y^2, z^2$, and $(yz)^2$, we 
see that $\pi_1(M(s))$ has as factor group 
\begin{eqnarray*}
\langle x,y,z &\mid& x^2, y^2, z^2, (yz)^2, (zx)^p, (yx)^q, \\
&&(yx)^{(q-1)/2}(zy)(yx)^{(q+1)/2}(zx)^{(p-1)/2}
(yz)(zx)^{(p+1)/2} \rangle .
\end{eqnarray*}

Let $G_{ev}$ denote the subgroup consisting of words of even length.
By substituting $a=zx$, $b=xy$, $c=yz$ we have the closely related group
$$G'_{ev} = \langle a,b,c \mid c^2, a^p, b^q,abc,b^{(q+1)/2}cb^{(q-1)/2}a^{(p-1)/2}
ca^{(p+1)/2} \rangle .$$
That is, $G_{ev}$ is a quotient $F_{a,b,c}/K_0$ of the free group on $a,b,c$ where
$K_0$ is  the normal closure of the given relations in the free group on $x,y,z$.
On the other hand, $G'_{ev} = F_{a,b,c}/K_1$ where $K_1$ is the normal closure of the
relations in $F_{a,b,c}$. Evidently $K_1 \leq K_0$ and it will suffice to argue that $G'_{ev}$ 
is infinite.

Letting $\al = a^{(p-1)/2}$, $\beta = b^{(q-1)/2}$ we have
$$\langle \al, \beta \mid \al^p, \beta^q, (\al^2 \beta^2)^2, 
(\beta \al)^2 (\al \beta)^2 \rangle.$$
Finally adding the relator $(\al \beta)^2$, we arrive at
$$(2,p,q; 2) = \langle \al, \beta \mid \al^p, \beta^q, (\al \beta)^2,
(\al^2 \beta^2)^2 \rangle.$$
\end{proof}

\subsection{ $2(p+q) \pm 1$ surgeries }

We will show that $2(p+q) \pm 1$ surgery results in a manifold with infinite fundamental group provided $p$ and $q$ are sufficiently large. 

\begin{lemma}
\label{lem2pqm1}%
Let $p,q$ be odd with $5 \leq p \leq q$ and $(p,q) \neq (5,5)$.
Let $M$ be the complement of the $(-2,p,q)$ pretzel knot.
Then $\pi_1(M(2(p+q)-1))$ is not finite.
\end{lemma}

\begin{proof}
As in the previous subsection, $\pi_1(M(2(p+q)-1))$  is 
\begin{eqnarray*} 
\langle x,y,z &\mid& 
(yz^{-1})^{-1}y(yz^{-1}) = (yx)^{(1-q)/2}x(yx)^{(q-1)/2},\\
&&(yz^{-1})^{-1}z(yz^{-1}) = (zx)^{(p+1)/2}x(zx)^{-(p+1)/2}, x^{2(p+q)-1}l
\rangle .
\end{eqnarray*}
We will write the group in terms of the generators $a = zx$ and $b = yx$.
Then $yz^{-1}$ = $ba^{-1}$. It will be convenient to use the expressions
$\alpha = a^{(p-1)/2}$ and $\beta = b^{(q-1)/2}$. 

The relator 
$x^{2(p+q)-1}l$ allows us to write $x$ in terms of $a$ and $b$:
$x = \beta a \beta \alpha b \alpha.$
Then, $y = bx^{-1} = b (\beta a \beta \alpha b \alpha)^{-1}$
while $z = ax^{-1} = a (\beta a \beta \alpha b \alpha)^{-1}$.
Using these substitutions, the first relation becomes the relator
$a \beta b^{-1} a \beta (\alpha b \alpha \beta)^2$
while the second yields
$\alpha b a^{-1} \alpha b (\alpha \beta a \beta)^2$.
Thus, $\pi_1(M(2(p+q)-1))$
 can be written
$$ \langle a, b, \alpha, \beta \mid
\alpha = a^{(p-1)/2}, \beta = b^{(q-1)/2},
a \beta b^{-1} a \beta (\alpha b \alpha \beta)^2,
\alpha b a^{-1} \alpha b (\alpha \beta a \beta)^2 \rangle .$$

As we shall see, by adding the relators
$a^p$, $b^q$, and $(\alpha \beta^{-1})^2$
we obtain $G^{5,p,q}$ (see~\cite{c}) 
as a factor group. Since this group is
infinite~\cite{ej}, we deduce that $\pi_1(M(2(p+q)-1))$ is, likewise, infinite.

Note that $a^p = a \alpha^2$ so that $a = \alpha^{-2}$. Similarly
$b = \beta^{-2}$. The relator $(\alpha \beta^{-1})^2$ gives
$\alpha \beta^{-1} = \beta \alpha^{-1}$ and $\alpha^{-1} \beta = \beta^{-1} \alpha$. Then, the relators of $\pi_1(M(2(p+q)-1))$
become $(\beta^{-1} (\alpha^2 \beta^{-2})^2)^2$ and
$((\alpha^2 \beta^{-2})^2 \alpha)^2$. Thus, we can write the factor
group as 
$$G = \langle \alpha, \beta \mid
\alpha^p, \beta^q, (\alpha \beta^{-1})^2,
(\beta^{-1} (\alpha^2 \beta^{-2})^2)^2,
((\alpha^2 \beta^{-2})^2 \alpha)^2 \rangle. $$
Replacing $\alpha$ by $A$, $\beta$ by $B^{-1}$, and introducing the generator $C = (A^2B^2)^2$, we can rewrite $G$ as
$$G = \langle A,B,C \mid
A^p, B^q, (AB)^2,
(BC)^2, (CA)^2, C = (A^2B^2)^2 \rangle.$$

We next show that $C^5$ is also a relator in this group. Since
$(BC)^2$ is a relator, so too is $(C^{-1}B^{-1})^2$. But
$$C^{-1}B^{-1} = (B^{-2}A^{-2})^2B^{-1} = B^{-1}(B^{-1}A^{-2}B^{-1})^2
= B^{-1}(AB^2A)^2.$$ It follows that
$(AB^2A)^{-2} = B^{-1}(AB^2A)^2B^{-1}$.

On the other hand, since 
$$1 = (CA)^2 = ((A^2B^2)^2A)^2 = (A(AB^2A)^2)^2$$
we have $(AB^2A)^{-2} = A(AB^2A)^2A$. Equating these two expressions
for $(AB^2A)^{-2}$ gives the desired result, $1 = (A^2B^2)^{10} = C^5$.

Note that $(ABC)^2$ is also a consequence of the relators we already know.
Indeed,
\begin{eqnarray*}
(ABC)^2 & =&  ABCABC \\
& = & B^{-1}A^{-1}A^{-1}C^{-1}C^{-1}B^{-1} \\
& = & B^{-1}A^{-2}C^{-2}B^{-1} \\
& = & B^{-1}A^{-2}(A^2B^2)^{-4}B^{-1} \\
& = & B (A^2B^2)^{-5}B^{-1} = 1.
\end{eqnarray*}

Thus, we will not change the group by adding the relators $C^5$ and
$(ABC)^2$:
$$G = \langle A,B,C \mid
A^p, B^q, C^5, (AB)^2,
(BC)^2, (CA)^2, (ABC)^2, C = (A^2B^2)^2 \rangle.$$
Finally, we note that $C = (A^2B^2)^2$ is now a consequence of the other relators. Indeed, as above, $(ABC)^2$ implies $C^{-2} = A^2B^2$
so that $C = C^{-4} = (A^2B^2)^2$. Thus, $G$ is the group
$G^{5,p,q} = \langle A,B,C \mid
A^p, B^q, C^5, (AB)^2,
(BC)^2, (CA)^2, (ABC)^2 \rangle$
defined by Coxeter~\cite{c}. Edjvet and Juh\'asz~\cite{ej} have recently shown that this group is infinite when $5 \leq p \leq q$, except in
the case where $p = q = 5$. Since $\pi_1(M(2(p+q)-1))$ projects onto
$G^{5,p,q}$, it is also infinite.

\end{proof}

\begin{rem}\label{rem:surjection}
We can use this argument to show $s = 2(p+q) - k$ surgery is infinite provided $k \equiv 1 \bmod 5$. The idea is to add the relator $x^{k-1}$ to form a factor group of $\pi_1(M(2(p+q)-k))$. 
Note that $l$ is of the form $x^{-2(p+q)}\tilde{l}$. 
After adding $x^{k-1}$, 
the relator $x^sl$ becomes $x^{-1}\tilde{l}$ as before and we arrive at the same group $G^{5,p,q}$ with an additional relator that corresponds to  $x^{k-1}$. The relation $x = \tilde{l}$ means that we can equally well think 
of adding the relator $\tilde{l}^{k-1}$. However, when we follow 
this through to the group $G^{5,p,q}$, we see that this amounts to 
adding the relator $C^{k-1}$. If $5 \mid (k-1)$ this relator is already
satisfied in $G^{5,p,q}$. Thus, for such $k$, $\pi_1(M(2(p+q)-k))$
projects onto $G^{5,p,q}$ and is infinite when
$5 \leq p \leq q$ and $(p,q) \neq (5,5)$.
For example, this shows that
$18$ surgery on the $(-2,5,7)$ knot and $22$ 
surgery on $(-2,5,9)$ are not finite.
\end{rem}

\begin{lemma}
\label{lem2pq1}%
Let $p,q$ be odd with $7 \leq p \leq q$. If $p = 7$, let $q \geq 21$ and 
suppose that $(p,q) \neq (9,9)$.
Let $M$ be the complement of the $(-2,p,q)$ pretzel knot.
Then $\pi_1(M(2(p+q)+1))$ is not finite.
\end{lemma}

\begin{proof}
As in the proof of Lemma~\ref{lem2pqm1},
$\pi_1(M(2(p+q)+1))$  is
$$ \langle a, b, \alpha, \beta \mid
\alpha = a^{(p-1)/2}, \beta = b^{(q-1)/2},
\beta \alpha b \alpha \beta (a \beta)^2 (\alpha b)^2,
\alpha \beta a \beta \alpha (a \beta)^2 (\alpha b)^2
\rangle . $$
Adding the relators $a^p$, $b^q$, $(\alpha \beta^{-1})^2$ yields
the group $G^{3,p,q}$. Edjvet and Juh\'asz~\cite{ej} have shown that this group is infinite under the given conditions on $p$ and $q$.
\end{proof}

\begin{rem} Edjvet and Juh\'asz determine finiteness of the groups $G^{m,p,q}$ except in the case $(m,n,p) \in \{(3,8,13), (3,7,19) \}$.
The argument above suggests that techniques of knot theory may be of use
in resolving these outstanding cases.
In particular, we've shown that $G^{3,7,19}$ is related to $2(7+19) + 1$ surgery of the $(-2,7,19)$ pretzel knot. We will argue
below that this surgery is not finite. Can this 
be used to resolve the open question of the finiteness of $G^{3,7,19}$?
\end{rem}

\section{\label{sec255}%
The $(-2,5,5)$ pretzel knot}

In this section, $K$ will denote the $(-2,5,5)$ pretzel knot and
we will prove
\begin{thm}
\label{thm255}
The $(-2,5,5)$ pretzel knot admits no non-trivial finite surgery.
\end{thm}

As a first step, we show that all but one of the boundary slopes of $K$ is detected.
By the Hatcher-Oertel~\cite{ho} algorithm, $K$ has boundary slopes $0$, $14$, $15$, $20$, and $22$. 
In the first subsection, using the method of~\cite{k}, we see that slopes $14$ and $15$ are detected
and  we show how that method can be extended to prove that the boundary slopes $20$ and $22$ are also detected and, 
moreover, there are at least two ideal points for $20$.

Having shown that all boundary slopes except $0$ are detected, we have a good idea of the Culler-Shalen norm of $K$. 
In the second subsection, we apply this knowledge toward a proof of Theorem~\ref{thm255}.

\subsection{Ideal points of the complement of the $(-2,5,5)$ pretzel knot}
We now explain how the technique of~\cite{k} (with which we assume familiarity)  
can be used to detect the $14$ and $15$ slopes.
By additional calculation we also show that each of the boundary slopes $20$ and $22$ are detected 
with slope $20$ having at least two ideal points.

SnapPea~\cite{weeks} gives an ideal triangulation of $S^3\setminus K$ with
$7$ ideal tetrahedra. 
Let $z_1, \dots , z_7$ be the complex parameters of the ideal tetrahedra.
We define $z'_k=\frac{1}{1-z_k}$ and $z''_k=1-\frac{1}{z_k}$.
The gluing equations associated to this ideal triangulation are written in the form
\[
\begin{split}
z_1 z'_2 z_3 z_4 z_5  &= z_1 (1-z_2)^{-1} z_3 z_4 z_5 =1, \\
z'_1 z''_1 z_2 z'_4 z'_5 z''_5 (z'_6)^2 z''_6 &= -(z_1)^{-1} z_2 (1-z_4)^{-1} (z_5)^{-1} (z_6)^{-1} (1-z_6)^{-1}=1, \\
z''_1 z_2 z''_2 z'_3 z''_3 z'_4 (z'_7)^2 z''_7 &=(z_1)^{-1}(1-z_1)(1-z_2) (z_3)^{-1}(1-z_4)^{-1}(z_7)^{-1} (1-z_7)^{-1}=1, \\
z_1 z'_2 z_6 z_7 &=z_1 (1-z_2)^{-1} z_6 z_7=1, \\
z_3 z''_4 z'_5 z_7 &=-z_3 (z_4)^{-1} (1-z_4) (1-z_5)^{-1} z_7=1, \\
z'_3 z''_4 z_5 z_6 &=-(1-z_3)^{-1} (z_4)^{-1} (1-z_4) z_5 z_6=1, \\
z'_1 z''_2 z''_3 z_4 z''_5 z''_6 z''_7 &=1 \\
\end{split}
\]
and the derivatives of the holonomies of the meridian and longitude are given by
\begin{equation}
\label{eq:ML}
\begin{split}
M & = z''_1 z'_4 (z_5)^{-1} z'_6, \\
L & = (z_1)^{-1} z'_1 (z''_1)^{-19} z_2 z'_2 (z'_3)^{-1} (z''_3)^{-1} 
(z'_4)^{-19} z''_4 (z_5)^{18} (z''_5)^{-1} (z'_6)^{-19}.
\end{split}
\end{equation}
Because the product of all the gluing equations is equal to $1$, we can omit the last equation.
Each gluing equation can be written in the form 
\[
\prod_{k=1}^7 (z_k)^{r'_{j,k}} (1-z_k)^{r''_{j,k}} = \pm 1 \quad (j = 1, \dots , 6)
\]
for some integers $r'_{j,k}$ and $r''_{j,k}$.
We denote by $\mathcal{D}$ the affine algebraic set defined by the gluing equations in 
$(\mathbb{C}-\{0,1\})^7$ and call it a {\em deformation variety}.
For a given point of $\mathcal{D}$, we can construct a $\PSLC$-representation by using the developing map.
It is known that this construction defines an algebraic map from $\mathcal{D}$ to the $\PSLC$-character variety.

As a sequence of points on $\mathcal{D}$ approaches an ideal point,
some $z_k$ goes to $0$, $1$, or $\infty$.
Let $I=\{i_1, \dots, i_7\}$ where $i_k$ is $1$, $0$, or $\infty$.
The vector $I$ represents a type of degeneration of the ideal tetrahedra.
Let 
\[
r(I)_{j,k}=
\left\{ \begin{array}{ll}
r''_{j,k}  & \textrm{if $i_k = 1$}\\
r'_{j,k}  & \textrm{if $i_k = 0$ }\\
-r'_{j,k}-r''_{j,k}  & \textrm{if $i_k = \infty$}
\end{array} \right. (k = 0, 1 \dots ,6) 
\]
and
\[
R(I)_{k}= \mathrm{det} 
\begin{pmatrix} 
r(I)_{1,1} & \hdots & \widehat{r(I)_{1,k}}   & \hdots & r(I)_{1,n} \\
\vdots   &        & \vdots     &        & \vdots \\
r(I)_{n-1,1} & \hdots & \widehat{r(I)_{n-1,k}}   & \hdots       & r(I)_{n-1,n}
\end{pmatrix}
\]
where the hat means removing the column.
We define $d(I)=(R(I)_0, -R(I)_1, \dots , R(I)_6)$. 
In~\cite{k}, it is shown that if all the coefficients of $d(I)$ are positive (or negative), 
there is a corresponding ideal point of $\mathcal{D}$.
The corresponding valuation $v$ satisfies $v(z_k)=d_k$.
So we can easily compute the values $v(M)$ and $v(L)$.
If one of $v(M)$ or $v(L)$ is non-zero, the ideal point of $\mathcal{D}$ gives an ideal point of  the 
$\PSLC$-character variety and its boundary slope is $-v(L)/v(M)$.
Because every $\PSLC$-representation of the fundamental group of a knot complement
lifts to a $\SLC$-representation,
we obtain an ideal point of the $\SLC$-character variety.

By computing $d(I)$ for all $I$, we obtain $6$ ideal points of $\mathcal{D}$ satisfying the condition.
The computation is shown in Figure \ref{fig:idealpoints}.
As in the figure, we conclude that the slopes $14$ and $15$ are detected.

\begin{figure}
\begin{tabular}[]{c|c|c|c}
$I$ & $d(I)$ & $(v(M),v(L))$ & boundary slopes \\ \hline
$(0,1,\infty,0,0,\infty,0)$ & $-(1,2,1,1,1,1,2)$ & $(1,-14)$ & $14$ \\ 
$(\infty,\infty,0,0,\infty,0,\infty)$ & $(2,1,1,1,1,2,1)$ & $(-1,14)$ & $14$ \\ 
$(\infty,1,\infty,0,0,\infty,0)$ & $-(1,4,2,4,3,1,6)$ & $(2,-30)$ & $15$ \\ 
$(0,1,\infty,\infty,0,\infty,0)$ & $-(4,3,1,1,1,2,1)$ & $(2,-30)$ & $15$ \\ 
$(\infty,1,0,0,\infty,0,\infty)$ & $(4,1,3,4,2,6,1)$ & $(-2,30)$ & $15$ \\ 
$(\infty,\infty,0,\infty,\infty,0,\infty)$ & $(3,4,1,1,1,1,2)$ & $(-2,30)$ & $15$ \\ 
\end{tabular}
  \caption{Ideal points of $\mathcal{D}$ detected by the method of~\cite{k}.}
  \label{fig:idealpoints}
\end{figure}

\subsubsection{Slope 20 (toroidal)}
If some ideal tetrahedron does not converge to $1$, $0$, or $\infty$,
we cannot apply the method of~\cite{k} directly.
Instead, we make a careful analysis of the
non-degenerate ideal tetrahedra.
When $(z_1,z_2,z_3,z_4,z_5,z_6,z_7)$ goes to $(*,0,*,*,*,1,*)$, there are corresponding ideal points.
(The $*$ means that the ideal tetrahedron does not become degenerate.)
We change the coordinate system of the degenerate ideal tetrahedra by  setting
\[
z_2=bt, z_6=1-ft.
\]
When $t=0$, the points correspond to ideal points.
There are two solutions:
\[
(z_1,z_3,z_4,z_5,z_7)=(-1, \frac{3\pm\sqrt{-3}}{2}, \frac{3\pm\sqrt{-3}}{6}, \frac{-1\pm\sqrt{-3}}{2}, -1).
\]
At the ideal point, we have $(v(z_2),v(1-z_2))=(1,0)$, $(v(z_6),v(1-z_6))=(0,1)$ 
and $(v(z_k),v(1-z_k))=(0,0)$ for $k \neq 2,6$. 
From the equations (\ref{eq:ML}), we have $v(M)=-1$ and $v(L)=20$.
So the corresponding boundary slope is $-v(L)/v(M)=20$.
We remark that the volume near the ideal points approaches 
$\pm  2.029883...= \pm 2\textrm{(the volume of the regular ideal tetrahedron)}$ respectively
in the sense of  section 8 of \cite{neumann-yang}.
So the representations near these ideal points are not conjugate. 
Therefore we conclude that these two ideal points of $\mathcal{D}$ give two different ideal points 
on the $\PSLC$-character variety.

We remark that 
when $(z_1,z_2,z_3,z_4,z_5,z_6,z_7)$ goes to $(1,*,*,*,*,*,1)$, there are also corresponding ideal points. 
While they also give a boundary slope of $20$, 
they appear to be equal to the above two ideal points on the $\PSLC$-character variety. 
(In general, the map from $\mathcal{D}$ to the $\PSLC$-character variety is two to one near an ideal point.) 

\subsubsection{Slope 22}
When $(z_1,z_2,z_3,z_4,z_5,z_6,z_7)$ goes to $(1,1,0,*,1,*,0)$, there is a corresponding ideal point.
We change the coordinate system for the degenerate ideal tetrahedra by setting
\[
z_1=1-at, z_2=1-bt, z_3=ct, z_5=1-et^2, z_7=gt.
\]
The solution is
\[
(z_4,z_6)=(\frac{1}{2},-1).
\]
The corresponding boundary slope is $22$.

When $(z_1,z_2,z_3,z_4,z_5,z_6,z_7)$ goes to $(\infty,0,1,*,0,0,*)$, there is also a corresponding ideal point.
This is equal to the above ideal point on the character variety.

\subsection{Finite surgeries of $(-2,5,5)$}

Using Hatcher and 
Oertel's algorithm~\cite{ho}, the boundary slopes of the $(-2,5,5)$ pretzel 
knot are $0$, $14$, $15$, $20$, and $22$. It follows that the Culler-Shalen norm is of the form
$$
\| \gamma_s \| = 2 {[} a_1 \Delta( s, 0)
+ a_{2} \Delta( s, 14)
+ a_{3} \Delta( s, 15)
+ a_{4} \Delta( s, 20)
+ a_{5} \Delta( s, 22) {]}
$$
with $a_i$ non-negative integers. Here $\gamma_s \in H_1(\partial M; \Z)$ is the class of a curve representing the slope $s \in \Q \cup \{\frac10\}$.

In the previous subsection we showed that
all boundary slopes other than $0$ are detected, and that $20$ is detected by two ideal points. It follows that the corresponding constants satisfy
$a_2, a_3, a_5 \geq 1$ and $a_4 \geq 2$. However, the following lemma
shows that $a_3 \geq 2$.

\begin{lemma}
\label{lemaeven}
Let $K$ be a hyperbolic knot in $S^3$ with Culler-Shalen norm
$\| \gamma \| = 2 \sum a_j \Delta(\gamma, \beta_j)$, the sum being taken over the finite set of boundary slopes, $B = \{ \beta_j \}$.
If $\beta_0 \in B$ is represented by a fraction $\frac{u}{v} \in \Q \cup \{ \frac10 \}$ with $u$ odd then the corresponding constant $a_0$ is even.
\end{lemma}

\begin{proof}
The A-polynomial of $K$ was defined in~\cite{ccgls} as an integral coefficient polynomial $A(l,m)$
in the variables $l$ and $m$. The authors also show that the 
monomials appearing in $A$ all have $m$ raised to an even power.
It follows that the Newton polygon of $A$ has vertices at lattice points $(x,y) \in \Z^2$ with $y$ even.

Boyer and 
Zhang~\cite{bz2} showed that the Newton polygon of $A$ is equivalent to
the Culler-Shalen norm.
Specifically, the vectors that connect consecutive vertices of the Newton polygon
have the form $a_j(v_j,u_j)$ where $\frac{u_j}{v_j}$ is the slope
corresponding to the boundary class  $\beta_j$. Under our hypotheses
then, the Newton polygon of $A$ would include 
vertices of the form $(x_0, y_0)$ and $(x_0 + a_0 v, y_0 + a_0 u)$. 
As both $y_0$ and $y_0 + a_0u$ are even, and $u$ is odd, we must have
that $a_0$ is even.
\end{proof}

Thus, using the results of the previous section, we can assume
$a_2, a_5 \geq 1$ and $a_3, a_4 \geq 2$.

\begin{proof}[Proof of Theorem~\ref{thm255}]
The algorithm of Hatcher and Oertel~\cite{ho} shows that the 
boundary slopes of the $(-2,5,5)$ pretzel knot are $0$, $14$, $15$, $20$, and $22$. By Lemma \ref{lemS}, the minimal norm of $K$ is $S=20$. Thus, using Lemmas~\ref{lem610}, \ref{lembz}, \ref{lemims}, and \ref{lemfinbdy}, a non-trivial finite surgery must lie in the set
$$\{ 13, 16, 19, 21, 23, \frac{31}{2}, \frac{39}{2}, \frac{41}{2}, \frac{43}{2}, \frac{45}{2} \}.$$
We will use the estimates of the $a_i$'s discussed above 
to show that none of these has small enough norm to be a finite slope. 

If $s \in \{\frac{39}{2}, \frac{41}{2}, \frac{43}{2}, \frac{45}{2} \}$,
then $\| \gamma_s \| \geq 2(11 a_2 + 9 a_3) \geq 58 > S+8$. Thus, $s$
is not a finite slope by Lemma~\ref{lembz}.
Similarly, if $s = \frac{31}{2}$, $\| \gamma_s \| \geq 2(3 a_2 + a_3 + 9 a_4 + 13 a_5) \geq 72 > S+8$.

If $s = 21$ or  $23$, 
then $\| \gamma_s \| \geq 2(7 a_2 + 6 a_3) \geq 38 > S+8$, and these
slopes are not finite. Also,
$\| \gamma_{13} \| \geq 2( a_2 + 2 a_3 + 7 a_4 + 9 a_5) \geq 56 > S+8$,
and $\| \gamma_{19} \| \geq 2( 5 a_2 + 4 a_3 + a_4 + 3 a_5) \geq 36 > S+8$,
so that $13$ and $19$ are likewise not finite.

This leaves only the slope $16$ as a candidate for finite surgery.
Since $\frac10$ is not a boundary slope, it follows from \cite[Corollary 1.1.4]{cgls} and Lemma~\ref{lemS} that $\| \frac10 \| = S = 20$.
We deduce that $\sum_{i=1}^5 a_i = 10$. Given our constraints on the
$a_i$'s, the smallest possible norm for slope $16$ is $44$, which occurs with the choice $(a_1, a_2, a_3, a_4, a_5) = (0, 1, 6, 2, 1)$. That is,
any other distribution of the $a_i$ consistent with $\sum_{i=1}^5 a_i = 10$, $a_2, a_5 \geq 1$, and $a_3, a_4 \geq 2$ will result in
$\| 16 \| \geq 44  > 40 = 2S$. Therefore, by Lemma~\ref{lembz}, $16$ is also not a finite slope.

Thus  we have shown that the only finite slope of the $(-2,5,5)$ pretzel knot is the trivial slope $\frac10$.
\end{proof} 

\section{\label{sec2579}%
The $(-2,5,7)$ and $(-2,5,9)$ pretzel knots}

In this section, we prove
\begin{thm}
\label{thm257}
The $(-2,5,7)$ pretzel knot admits no non-trivial finite surgery.
\end{thm}
and 
\begin{thm}
\label{thm259}
The $(-2,5,9)$ pretzel knot admits no non-trivial finite surgery.
\end{thm}

We prove these theorems using information about detected boundary slopes.
For Theorem~\ref{thm259}, it suffices to use the slopes
detected by the methods of~\cite{k}. For Theorem~\ref{thm257}, however, we
will need to also know something about the number of ideal points
for the boundary slope $24$.
For this, we will use Ohtsuki's method for calculating the number of ideal points~\cite{o1,o2}. 
In the first subsection we give a brief overview of Ohtsuki's approach and
find lower bounds for the number of ideal points for the boundary slope $24$ of the 
$(-2,5,7)$ pretzel knot.
We then give proofs of our theorems in the subsequent two subsections.

\subsection{Detecting ideal points using Ohtsuki's method}

In this subsection we show that
there are at least $8$ ideal points for the boundary slope $24$ of the $(-2,5,7)$ pretzel knot.
These ideal points are detected by using Ohtsuki's method as outlined in~\cite{o1,o2}. 
We first briefly recall Ohtsuki's method and then give the calculation of the number of ideal points.

Fix a diagram of a knot $K$ in $S^3$ and consider the Wirtinger presentation  of 
$\pi_1(M)=\pi_1(S^3\setminus N(K))$
with generators $x_1,\cdots,x_c$, where $c$ is the number of crossings in the diagram.
For $\rho\in\text{Hom}(\pi_1(M),SL(2,\C))$, we set $X_i=\rho(x_i)$ for $i=1,\cdots,c$.
Since all $X_i$'s are conjugate, their eigenvalues are the same, say $\lambda$ and $\frac{1}{\lambda}$.
We assume that $|\lambda|>1$ since we are only interested in the case where $\lambda$ diverges.
Let $\begin{pmatrix}x_i^+ \\ 1\end{pmatrix}$
and $\begin{pmatrix}x_i^- \\ 1\end{pmatrix}$ denote the eigenvectors of $X_i$ corresponding to
$\lambda$ and $\frac{1}{\lambda}$ respectively.
For each relation $x_k=x_ix_jx_i^{-1}$ in the Wirtinger presentation
we have the following two equations:
\[
\begin{split}
   &(x_i^\pm-x_j^-)(x_k^\pm-x_j^+)-\ve (x_i^\pm-x_j^+)(x_k^\pm-x_j^-)=0,
\end{split}
\]
where $\ve=\lambda^{-2}$. Since there are $c$ relations, we have $2c$ equations of this type.
We denote them by $R_1,\cdots,R_{2c}$.
Let $\hat R(M)$ denote the algebraic set in $\C^{2c+1}\cap\{|\ve|<1\}$ 
with coordinates $(x_1^+,x_1^-,\cdots,x_c^+,x_c^-,\ve)$ determined by $R_1=\cdots=R_{2c}=0$.
Using a M\"obius transformation on $\C$, we can fix three variables
in $(x_1^+,x_1^-,\cdots,x_c^+,x_c^-)$, which corresponds to taking a slice of $\hat R(M)$ with three hyperplanes.
We denote this slice by $s(M)$.
The character variety $\chi(M)$ is defined to be the set of characters
of the $SL(2,\C)$-representations of $\pi_1(M)$.
There is a canonical map $t:\hat R(M)\to\chi(M)\cap\{|\ve|<1\}$ defined by 
$t(X)=\text{trace}(X)$.
By restricting this map to $s(M)$, we have a surjection $t_s:s(M)\to\chi(M)\cap\{|\ve|<1\}$.

Below, we summarize the algorithm for finding ideal points of $\chi(M)$ 
in~\cite{o1,o2}.
\vspace{3mm}

\noindent
{\bf Step 1.}\;\,
Suppose that we wish to detect 
an essential surface $S$ with a given boundary slope.
Choose a set of loops $\ell_1,\cdots,\ell_m$ in $S^3\setminus S$
and make a tree with ends $x_i^\pm$ according to~\cite{o2}
(see for instance Section~2 in~\cite{o2}) in such a way that
each of $\ell_1,\cdots,\ell_m$ has a fixed point when acting on the tree.
Note that when we make the tree using Ohtsuki's argument in~\cite{o2},
we use all the equations $R_1=\cdots,R_{2c}=0$ 
to determine $\hat R(M)$, without omitting three of them. Instead, the three are removed after Step 3 below.
\vspace{3mm}

\noindent
{\bf Step 2.}\;\,
Fix three variables using the M\"obius transformation.
For the other variables, according to the tree obtained in Step~1,
we change the coordinates as follows:
if the tree suggests that $x_i^\bullet\equiv x_j^\bullet \mod\ve$
then we introduce
a new variable $\tilde x_j^\bullet$ in place of $x_j^\bullet$ such that 
$x_j^\bullet=x_i^\bullet+\tilde x_j^\bullet\ve$,
where~$\bullet$ represents~$+$ or~$-$ (and the choice can be made 
independently for $x_i^\bullet$ and $x_j^\bullet$).
Let $\tilde R_1,\cdots,\tilde R_{2c}$ denote
the equations of $R(M)$ in the new variables.
These determine an algebraic set in $\C^{2c+1}\cap\{|\ve|<1\}$ with new coordinates
and we denote it by $\tilde R(M)$.
The three fixed variables determine a slice of $\tilde R(M)$ and we denote it by $\tilde s(M)$.
Note that the above change of coordinates determines a canonical map $\phi:\tilde R(M)\to \hat R(M)$
and its restriction $\phi|_{\tilde s(M)}:\tilde s(M)\to s(M)$.
\vspace{3mm}

\noindent
{\bf Step 3.}\;\,
According to~\cite[Section 2.4]{o2}, we make a tunnel near some crossing.
Let $R_1, R_2, R_3$ and $R_4$  be the equations 
which have to be removed due to the tunnel.
Suppose that the generators $x_i$ (resp.\ $x_j$) of $\pi_1(M)$ is separated 
into $x_i$ and $x_i'$ (resp. $x_j$ and $x_j'$) by the tunnel as in~\cite[Fig.11]{o2}.
\vspace{3mm}

We will omit the equations $R_1$, $R_2$, and $R_3$.
\vspace{3mm}

\noindent
{\bf Step 4.}\;\,
Let $\hat R'(M)$ denote the algebraic set in $\C^{2c+1}\cap\{|\ve|<1\}$
determined by $R_4=\cdots=R_{2c}=0$
and $\tilde R'(M)$ denote the one determined by $\tilde R_4=\cdots=\tilde R_{2c}=0$.
The change of coordinates in Step~2 determines
a canonical map $\phi':\tilde R'(M)\to \hat R'(M)$.
The three fixed variables determine 
slices $s'(M)$ in $\hat R'(M)$ and $\tilde s'(M)$ in $\tilde R'(M)$.
\vspace{3mm}

\begin{prop}\label{mainlemma}
Let $z$ be an isolated point of the set $\tilde s'(M)\cap\{\ve=0\}$ and set $w=\phi'(z)$.
Suppose that the generators $x_i,\,x_j\in\pi_1(M)$ separated by the tunnel
satisfy $x_i^+\ne x_j^+$, $x_i^+\ne x_j^-$, 
$x_i^-\ne x_j^+$ and $x_i^-\ne x_j^-$ in a neighbourhood of $w$.
Then $w$ is an ideal point of $s(M)$.
\end{prop}

%
%

\begin{proof}
We first show that there is a curve in $\tilde s'(M)\setminus\{\ve=0\}$ 
which converges to $z$ as $\ve\to 0$.
The set $\tilde s'(M)\cap\{\ve=0\}$ is an algebraic set in $\C^{2c+1}$ determined by
the equations $R_4=\cdots=R_{2c}=0$, $\ve=0$ and three equations for making the slice.
Hence the number of equations is $2c+1$.
Since the solutions are isolated, we can say that this algebraic set is 
a locally complete intersection. Hence $\dim \tilde s'(M)=1$.
This curve cannot be locally contained in $\{\ve=0\}$ because if it were then locally 
we would have $\dim \tilde s'(M)=1$, which contradicts the assumption that
$z$ is isolated. Hence this curve satisfies the property claimed.

Now we prove the original assertion.
Using the hypotheses and~\cite[Lemma~2.5]{o2}, we have $X_i=X_i'$ and $X_j=X_j'$.
Thus we can recover the equations $\tilde R_1,\tilde R_2$, and $\tilde R_3$ 
at each point in a small neighbourhood of $z$.
This means that $\tilde s'(M)$ locally coincides with $\tilde s(M)$.
Hence the curve obtained in the first paragraph
can be regarded as a curve in $\tilde s(M)\setminus\{\ve=0\}$ which converges to $z$ as $\ve\to 0$.
By using the map $\phi|_{\tilde s(M)}:\tilde s(M)\to s(M)$,
we conclude that $w$ is an ideal point of $s(M)$.
\end{proof}

Since $\lambda$ diverges at the point $w$ in Proposition~\ref{mainlemma},
using the map $t_s:s(M)\to \chi(M)\cap\{|\ve|<1\}$,
we conclude that $w$ corresponds to an ideal point of $\chi(M)$.

\begin{rem}
The generators $x_i,\,x_j\in\pi_1(M)$ separated by the tunnel
both appear in one relation in the Wirtinger presentation. Hence,
if the tree in Step~1 is made from local pieces as shown in~\cite[Fig.~4]{o2},
i.e., if there is no degeneration, 
then we automatically have the condition 
$x_i^+\ne x_j^+$, $x_i^+\ne x_j^-$, $x_i^-\ne x_j^+$ and $x_i^-\ne x_j^-$.
\end{rem}

\begin{lemma}\label{cross_ratios}
If two solutions $\zeta$ and $\eta$ have different values in the complex ratio
\[
   \frac{x_i^\bullet-x_k^\bullet}{x_i^\bullet-x_\ell^\bullet}\left/\frac{x_j^\bullet-x_k^\bullet}{x_j^\bullet-x_\ell^\bullet}\right.,
\]
where $x_i^\bullet, x_j^\bullet, x_k^\bullet, x_\ell^\bullet$ are some ends
of generators $x_1,\cdots,x_c$,
then they correspond to different ideal points of $\chi(M)$.
\end{lemma}

\begin{proof}
Let $z$ denote the point in $s(M)$ corresponding to the solution $\zeta$ 
and $w$ denote the one corresponding to $\eta$.
Suppose that they correspond to the same ideal point of $\chi(M)$. Then
their $SL(2,\C)$-representations must be conjugate, i.e.,
there are neighbourhoods $z\in U\subset s(M)$ and $w\in V\subset s(M)$
such that there is a map $\varphi$ of conjugation satisfying
$\varphi(z)=w$ and $\varphi(U)=V$. This means that,
for each pair of points $u\in U$ and $v=\varphi(u)\in V$,
there exists an $SL(2,\C)$-matrix $X$ such that
$X\rho_u(a)X^{-1}=\rho_v(a)$ for every element $a$ in $\pi_1(M)$,
where $\rho_u(a)$ is the $SL(2,\C)$-representation of $a$ at $u\in U$
and $\rho_v(a)$ is the one of $a$ at $v=\varphi(u)\in V$.
As mentioned in the proof of Lemma~C.1 in~\cite{o2},
the M\"obius transformation by $X$ sends the ends of generators at $z$ to
those at $w$. Thus the complex ratios have to be the same for $\zeta$ and $\eta$.
\end{proof}

\begin{lemma}
There are at least $8$ ideal points for the boundary slope $24$ of the $(-2,5,7)$ pretzel knot.
\end{lemma}

\begin{proof}
Set the generators $a,b,c,d,s_1,\cdots,s_7,t_1,\cdots,t_9$ as shown 
on the left in Figure~\ref{fig2} with
identification $a=t_1$, $b=s_2$, $c=t_8$, $d=s_7$, $s_1=t_2^{-1}$ and $s_6=t_9^{-1}$.
\begin{figure}[htbp]
\centerline{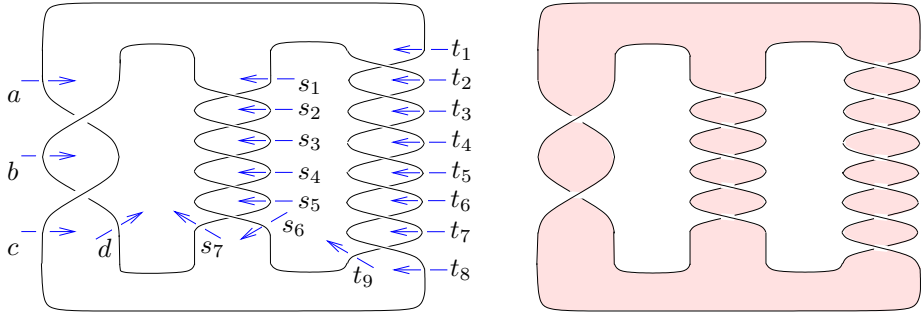}
   \caption{The $(-2,5,7)$ pretzel knot and an essential surface
with boundary slope $24$.\label{fig2}}
\end{figure}
We fix three parameters $t_1^-=0, t_1^+=1, t_3^+=t_2^++1$ and 
change the coordinate $t_2^+$ into $\zeta$ by setting $t_2^+=1+\zeta$.
We first consider Step~1.
Since the essential surface with boundary slope $24$ is as shown
on the right in Figure~\ref{fig2},
a tree consisting of subtrees of parasol type with the same origin
satisfies the necessary condition for corresponding to this surface. 
For the definition of a subtree of parasol type, see~\cite{o1}. 
We will find $8$ ideal points from such a tree.
The change of coordinates in Step~2 is done as
$b^-=a^++r_1\ve$, $c^-=b^++r_2\ve$, $d^-=c^++r_3\ve$,
$s_1^+=s_1^-+x_1$, $s_i^+=s_{i-1}^++x_i$ for $i=2,\cdots,7$, $s_i^-=s_{i-1}^++y_{i-1}\ve$ for $i=2,\cdots,7$,
$t_i^+=t_{i-1}^++p_i$ for $i=4,\cdots,9$, and $t_i^-=t_{i-1}^++q_{i-1}\ve$ for $i=2,\cdots,9$, 
where the $r_i$'s, $x_i$'s $y_i$'s, $p_i$'s, and $q_i$'s are non-zero.
Using the identification of tangles, we can set $a^\pm=t_1^\pm$, $b^+=s_2^+$, 
$s_1^-=\text{\tt solve}(b^--s_2^-,s_1^-)$, $c^+=t_8^+$,
$x_2=\text{\tt solve}(c^--t_8^-,x_2)$, $d^+=s_7^+$,
$x_3=\text{\tt solve}(d^--s_7^-,x_3)$,
$y_1=\text{\tt solve}(s_1^+-t_2^-,y_1)$,
$x_1=\text{\tt solve}(s_1^--t_2^+,x_1)$,
$y_6=\text{\tt solve}(s_6^+-t_9^-,y_6)$ and
$x_6=\text{\tt solve}(s_6^--t_9^+,x_6)$,
where {\tt solve} is the function which gives the value of the 
second entry that will make the first entry equal to zero.

The remaining variables are now
$\ve$, $\zeta$, $r_1, r_2, r_3$, $x_4, x_5, x_7$, $y_2,\cdots,y_5$,
$p_4,\cdots, p_9$ and $q_1,\cdots,q_8$, and hence $\tilde R(M)$ is an algebraic set in $\C^{26}$.
From the Wirtinger presentation, we have $28$ equations
$f_1^\pm$, $f_2^\pm$, $g_i^\pm$ for $i=1,\cdots,5$ and $h_i^\pm$ for $i=1,\cdots,7$
as shown in Figure~\ref{fig1}.
Let $F_1^\pm$, $F_2^\pm$, $G_i^\pm$ for $i=1,\cdots,5$ and $H_i^\pm$ for $i=1,\cdots,7$ denote
the leading coefficients of their expansions by $\ve$ after the change of coordinates in Step~2.
\begin{figure}[htbp]
\centerline{
\begin{picture}(0,0)%
\includegraphics{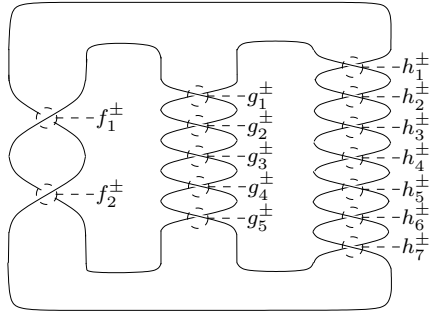}%
\end{picture}%
\setlength{\unitlength}{3158sp}%
\begingroup\makeatletter\ifx\SetFigFont\undefined%
\gdef\SetFigFont#1#2#3#4#5{%
  \reset@font\fontsize{#1}{#2pt}%
  \fontfamily{#3}\fontseries{#4}\fontshape{#5}%
  \selectfont}%
\fi\endgroup%
\begin{picture}(3173,2444)(1115,-5183)
\put(1882,-3705){\makebox(0,0)[lb]{\smash{{\SetFigFont{8}{9.6}{\rmdefault}{\mddefault}{\updefault}$f_1^\pm$}}}}
\put(1882,-4305){\makebox(0,0)[lb]{\smash{{\SetFigFont{8}{9.6}{\rmdefault}{\mddefault}{\updefault}$f_2^\pm$}}}}
\put(3069,-3533){\makebox(0,0)[lb]{\smash{{\SetFigFont{8}{9.6}{\rmdefault}{\mddefault}{\updefault}$g_1^\pm$}}}}
\put(3069,-3766){\makebox(0,0)[lb]{\smash{{\SetFigFont{8}{9.6}{\rmdefault}{\mddefault}{\updefault}$g_2^\pm$}}}}
\put(3069,-4006){\makebox(0,0)[lb]{\smash{{\SetFigFont{8}{9.6}{\rmdefault}{\mddefault}{\updefault}$g_3^\pm$}}}}
\put(3076,-4246){\makebox(0,0)[lb]{\smash{{\SetFigFont{8}{9.6}{\rmdefault}{\mddefault}{\updefault}$g_4^\pm$}}}}
\put(3069,-4493){\makebox(0,0)[lb]{\smash{{\SetFigFont{8}{9.6}{\rmdefault}{\mddefault}{\updefault}$g_5^\pm$}}}}
\put(4284,-3306){\makebox(0,0)[lb]{\smash{{\SetFigFont{8}{9.6}{\rmdefault}{\mddefault}{\updefault}$h_1^\pm$}}}}
\put(4284,-3539){\makebox(0,0)[lb]{\smash{{\SetFigFont{8}{9.6}{\rmdefault}{\mddefault}{\updefault}$h_2^\pm$}}}}
\put(4284,-3779){\makebox(0,0)[lb]{\smash{{\SetFigFont{8}{9.6}{\rmdefault}{\mddefault}{\updefault}$h_3^\pm$}}}}
\put(4284,-4019){\makebox(0,0)[lb]{\smash{{\SetFigFont{8}{9.6}{\rmdefault}{\mddefault}{\updefault}$h_4^\pm$}}}}
\put(4284,-4266){\makebox(0,0)[lb]{\smash{{\SetFigFont{8}{9.6}{\rmdefault}{\mddefault}{\updefault}$h_5^\pm$}}}}
\put(4284,-4486){\makebox(0,0)[lb]{\smash{{\SetFigFont{8}{9.6}{\rmdefault}{\mddefault}{\updefault}$h_6^\pm$}}}}
\put(4284,-4719){\makebox(0,0)[lb]{\smash{{\SetFigFont{8}{9.6}{\rmdefault}{\mddefault}{\updefault}$h_7^\pm$}}}}
\end{picture}%
}
   \caption{Positions of crossings corresponding to
the relations in the Wirtinger presentation.\label{fig1}}
\end{figure}
We omit the three equations $G_1^+, G_1^-$ and $G_2^-$ and 
solve the other equations in the following order:
$(H_1^+,q_1)$, 
$(H_1^-,q_2)$, 
$(H_2^-,q_3)$, 
$(H_3^-,q_4)$, 
$(H_4^-,q_5)$, 
$(H_5^-,q_6)$,
$(H_6^-,q_7)$,
$(H_7^-,q_8)$,
$(F_2^+,r_2)$, 
$(G_3^+,y_3)$, 
$(G_4^+,y_4)$, 
$(G_4^-,y_5)$,
$(G_5^+,x_4)$,
$(G_5^-,r_3)$,
$(H_3^+,p_5)$, 
$(H_2^+,p_4)$, 
$(G_3^-,x_7)$,
$(H_6^+,p_8)$,
$(H_4^+,p_6)$,
$(H_5^+,p_7)$,
$(H_7^+,p_9)$,
$(F_2^-,x_5)$,
$(F_1^+,r_1)$,
$(G_2^+,y_2)$,
where we again solve by finding the value of the second variable that will make the first zero.
Then checking the equation $F_1^-=0$ we have
\begin{eqnarray*}
0 & = & (\zeta^8+8\zeta^7+21\zeta^6+14\zeta^5-19\zeta^4-20\zeta^3+5\zeta^2+2\zeta+1) \\
& & \mbox{ } \times (\zeta^8+8\zeta^7+21\zeta^6+14\zeta^5-19\zeta^4-20\zeta^3+7\zeta^2+6\zeta-1). 
\end{eqnarray*}
By using {\tt resultant}, we can verify that the
$r_i$'s, $x_i$'s $y_i$'s, $p_i$'s, and $q_i$'s are non-zero for all solutions of this equation.
It is also easy to check that the $16$ solutions of this equation have 
$8$ different values in the cross ratio
\[
   \frac{t_1^+-t_2^+}{t_1^+-t_1^-}\left/\frac{t_3^+-t_2^+}{t_3^+-t_1^-}\right..
\]
Hence there are at least $8$ ideal points by Lemma~\ref{cross_ratios}.
\end{proof}

\begin{rem}
The methods of the above lemma can also be used to show that 
there are at least two ideal points for the boundary slope $20$ of the $(-2,5,5)$ pretzel knot.
In this case, we can find two other ideal points coming from an $r$-curve by the following reason:
Let $M$ denote the $(-2,5,5)$ pretzel knot complement and 
$M_{2,5}$ denote the $(2,5)$ torus knot complement.
There is a branched double covering map $M\to M_{2,5}$ such that
the essential surface $S$ with boundary slope $20$ 
is mapped to the essential M\"obius band in $M_{2,5}$ with boundary slope $10$. 
This means that $M$ has an $r$-curve with $r=20$.
Since it is shown in~\cite[Section~5.1]{o2} that
the boundary slope $10$ of $M_{2,5}$ has two ideal points,
we can conclude that this $r$-curve also has two ideal points.
The existence of these two additional ideal points can also be checked 
by using Ohtsuki's method for $M$ with a degenerate tree in Step~1.
Thus we can conclude that  the boundary slope~$20$ of the $(-2,5,5)$ pretzel knot
has at least four ideal points in total.
\end{rem}

\subsection{Finite surgeries on $(-2,5,7)$}

In this subsection, we prove Theorem~\ref{thm257}.

\begin{proof}[Proof of Theorem~\ref{thm257}]
Using~\cite{ho},
the norm of $(-2,5,7)$ is of the form
$$ \| \gamma \| = 2 {[} a_1 \Delta(\gamma, 0) +
a_{2} \Delta(\gamma, 14) + a_{3} \Delta(\gamma, 15) + 
a_{4} \Delta(\gamma, \frac{37}{2}) +
a_{5} \Delta(\gamma, 24) + a_{6} \Delta(\gamma, 26) {]} $$
where the $a_i$ are non-negative integers. By Lemma~\ref{lemS}, the 
minimal total norm is $S = 34$.
Using~\cite{k}, slopes $14$, $15$, and $\frac{37}{2}$ are detected. This implies $a_2 \geq 1$ and, by Lemma~\ref{lemaeven},
$a_3$ and $a_4$ must be even, so $a_3, a_4 \geq 2$.
As argued in the previous subsection, $24$ is detected by $8$ ideal points, 
so $a_5 \geq 8$.

Using Lemmas~\ref{lem610}, \ref{lembz}, \ref{lemims}, and \ref{lemfinbdy} a non-trivial finite surgery must lie in the set
$$\{ 16, 17, 18, 19, 20, 23, 25, 27, \frac{47}{2}, \frac{49}{2}, \frac{51}{2}, \frac{53}{2} \}.$$

If $s \in  \{ \frac{47}{2}, \frac{49}{2}, \frac{51}{2}, \frac{53}{2} \}$,
then $\| \gamma_s \| \geq 2 (19 a_2  + 17 a_3 + 20 a_4) \geq 186 > 42 = S+8$
and such a slope is not finite.
If $s \in \{23, 25, 27 \}$,
then $\| \gamma_s \| \geq 2 (9 a_2 + 8 a_3 + 9 a_4) \geq 86 > 42 = S+8$
and such a slope is not finite. 

For $s \in \{17, 19\}$, 
$\| \gamma_s \| \geq 2 (5 a_5) \geq  80 > 42 = S+8$ and these are not finite slopes. 
If $s \in \{16, 18 \}$,
$\| \gamma_s \| \geq 2 (6 a_5) \geq  96 > 68 = 2S$, so these are also
not finite.
Finally, 
$\|  \gamma_{20} \|  \geq 2 (6a_2 + 5 a_3 + 3 a_4  + 4 a_5) \geq 108 > 2S$.

Thus, the only finite slope of the $(-2,5,7)$ knot is the 
trival slope $\frac10$.
 \end{proof}






\subsection{Finite surgeries on $(-2,5,9)$}

In this subsection, we prove Theorem~\ref{thm259}.

\begin{proof}[Proof of Theorem ~\ref{thm259}]
As the argument is quite similar to that used to prove
Theorems~\ref{thm255} and \ref{thm257}, we will omit some details.
Using~\cite{ho} the norm is of the form
$$ \| \gamma \| = 2 {[} a_1 \Delta(\gamma, 0) +
a_{2} \Delta(\gamma, 14) + a_{3} \Delta(\gamma, 15) + 
a_{4} \Delta(\gamma, \frac{67}{3}) +
a_{5} \Delta(\gamma, 28) + a_{6} \Delta(\gamma, 30) {]}, $$
and $S = 48$ by Lemma~\ref{lemS}.
By~\cite{k}, slopes $14$, $15$, and $\frac{67}{3}$ are detected, so that $a_2 \geq 1$, and, by Lemma~\ref{lemaeven}, $a_3, a_4 \geq 2$.

Using the lemmas of Section~\ref{seclem}, the candidates for a finite slope are 
$$\{21, 22, 23, 24, 27, 29, 31, \frac{55}{2}, \frac{57}{2}, \frac{59}{2},  \frac{61}{2} \} .$$
The half integral surgeries will have norm at least 278 and slopes $27$, $29$, $31$ will exceed 130. Also, $\| 23 \| \geq 58$. 
These are all more than $S+8$, so none of these slopes are finite. 

For $21$ we must consider which distribution of the $a_i$'s will give the
least value for $\| 21 \|$. Since $\frac10$ is not a boundary slope, 
we have $24 = S/2 = \| \frac10 \|/2 = a_1 + a_2 + a_3 + 3 a_4 + a_5 + a_6$.
Thus, $\| 21 \|$ is minimised by the choice 
$\vec{a} = (a_1, a_2, a_3, a_4, a_5, a_6) = (0,2,4,6,0,0)$ which results in 
$\|21 \| =  124 > S+8$, so $21$ is not finite.
Similarly, we can minimise $\| 24 \|$ by choosing $\vec{a} = (0,1,2,6,3,0)$. Then $\|24 \| = 140 > 2S$, so $24$ is not finite.

Finally, as in Remark \ref{rem:surjection},
$\pi_1(M(22))$ surjects onto $G^{5,5,9}$ and is therefore not finite.
\end{proof}

\section{\label{sec6thm}%
The 6-Theorem}

In this section, let $K$ be a $(-2,p,q)$ pretzel knot with $p,q$ odd and either $7 \leq p \leq q$ or else $p=5$ and $q \geq 11$.
We will prove that $K$ admits no non-trivial finite surgeries. 

The argument proceeds in two steps: first, we reduce the candidate slopes for finite surgeries to a short list, and then we treat the short list. In the first subsection, 
we use the 6-theorem of Agol~\cite{a} and  Lackenby~\cite{l} to show that,
in case $7 \leq p \leq q$,   the candidates for a finite slope of $K$ are the trivial slope 
$\frac10$ and the integral slopes $2(p+q)+k$ with $k = -1, 0, 1, 2$.
Then, by Lemma~\ref{lemfinbdy}, the boundary slopes~\cite{ho}
$2(p+q)$ and $2(p+q)+2$ are not finite slopes. Lemma~\ref{lem2pqm1}
shows that $2(p+q)-1$ is not a finite slope while
Lemma~\ref{lem2pq1} says the same of the slope $2(p+q)+1$ except 
for eight cases with $p = 7$ or $9$. We will use Culler-Shalen norm arguments to address these cases in Section ~\ref{secpseven}.

In the second case, when $p = 5$ and $q \geq 11$, the 6-theorem implies that the only candidates for a finite slope of $K$ are $\frac10$ and the integral slopes $2(p+q)+k$ with $k=-2,-1,0,1,2,3$.
Lemma~\ref{lemfinbdy} again eliminates
$2q+10$ and $2q+12$, while 
the slopes $2q+8$ and $2q+9$ are also ruled out by Lemmas~\ref{lemeven} and \ref{lem2pqm1}. In Section~\ref{secpfive}, we will use the Culler-Shalen norm to rule out the slopes $2q+11$ and $2q+13$.

\subsection{Applying the $6$-theorem}
In  this section we shall study exceptional surgeries of the $(-2,p,q)$ pretzel knot complement.
Recall that the surgery on a knot $K$ along slope $s$ is called \emph{exceptional} if $K(s)$ is reducible, toroidal, or Seifert fibered, or if $\pi_1(K(s))$ is finite or not word-hyperbolic.  The following theorem was shown by Agol  \cite{a} and Lackenby \cite{l}.

\begin{thm}\label{6theorem}
Let $M$ be a compact orientable 3-manifold with interior 
having a complete, finite volume hyperbolic structure. 
Let $s_1, \dots , s_n$ be slopes on $\partial M$, with one $s_i$ on each component of $\partial M$. 
Suppose that there
is a horoball neighbourhood $N$ of the cusps of $M \setminus \partial M$ on which each $s_i$
has length more than $6$. 
Then, the manifold obtained by Dehn surgery along
$s_1, \dots , s_n$ is not exceptional: it is irreducible, atoroidal and not Seifert fibred, and has infinite
word hyperbolic fundamental group.
\end{thm}

By the geometrization theorem, all non-exceptional surgeries yield hyperbolic manifolds.
However, for our purposes we only need the conclusion that the filled manifold has infinite fundamental group. We shall apply Theorem  \ref{6theorem} to prove the following two propositions.

\begin{prop}\label{propk1}
For $p, q\geq 7$ and odd, the $(-2,p,q)$ pretzel knot has at most five exceptional surgeries.
If $s$ is an exceptional surgery, it is one of the following:
 $\frac10$, $2(p+q)-1$, $2(p+q)$, $2(p+q)+1$, or $2(p+q)+2$.
\end{prop}

\begin{prop}\label{propk2}
For $q \geq 11$ and odd, the $(-2,5,q)$ pretzel knot has at most seven exceptional surgeries. If $s$ is an exceptional surgery, it is one of the following:
$\frac10$, $2(5+q)+k$ $(k=-2,-1,0,1,2,3)$. 
\end{prop}

\begin{rem}
As shown in \cite{m1}, $2(p+q)$-surgery is a toroidal surgery and
the trivial surgery $\frac10$ is also exceptional. However, in general,
the other surgeries listed in the propositions may be hyperbolic.

\end{rem}

\begin{proof}[Proof of Proposition~\ref{propk1}]
Let $L$ be the link consisting of the $(-2,1,1)$ pretzel knot and 
two trivial link components encircling the $p$-twist and the $q$-twist. (See Figure \ref{fig:thelink}.)
By $-1/k$-surgery and $-1/l$-surgery along the two trivial link components 
we get the $(-2,1+2k,1+2l)$ pretzel knot complement in $S^3$.  
We will construct an ideal triangulation of $S^3\setminus L$ and find its complete hyperbolic structure.
Then we will study cusp shapes of the link complement
and apply the 6-theorem to that 3-cusped manifold.

\begin{figure}
\begin{center}
\includegraphics[width=80pt,clip]{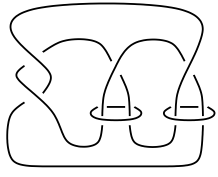}
\caption{By $-1/k$-surgery and $-1/l$-surgery along the trivial link components we get
the $(-2,1+2k,1+2l)$ pretzel knot.}
\label{fig:thelink}
\end{center}
\end{figure}

\begin{figure}[hb]
\begin{center}
\includegraphics[width=60pt,clip]{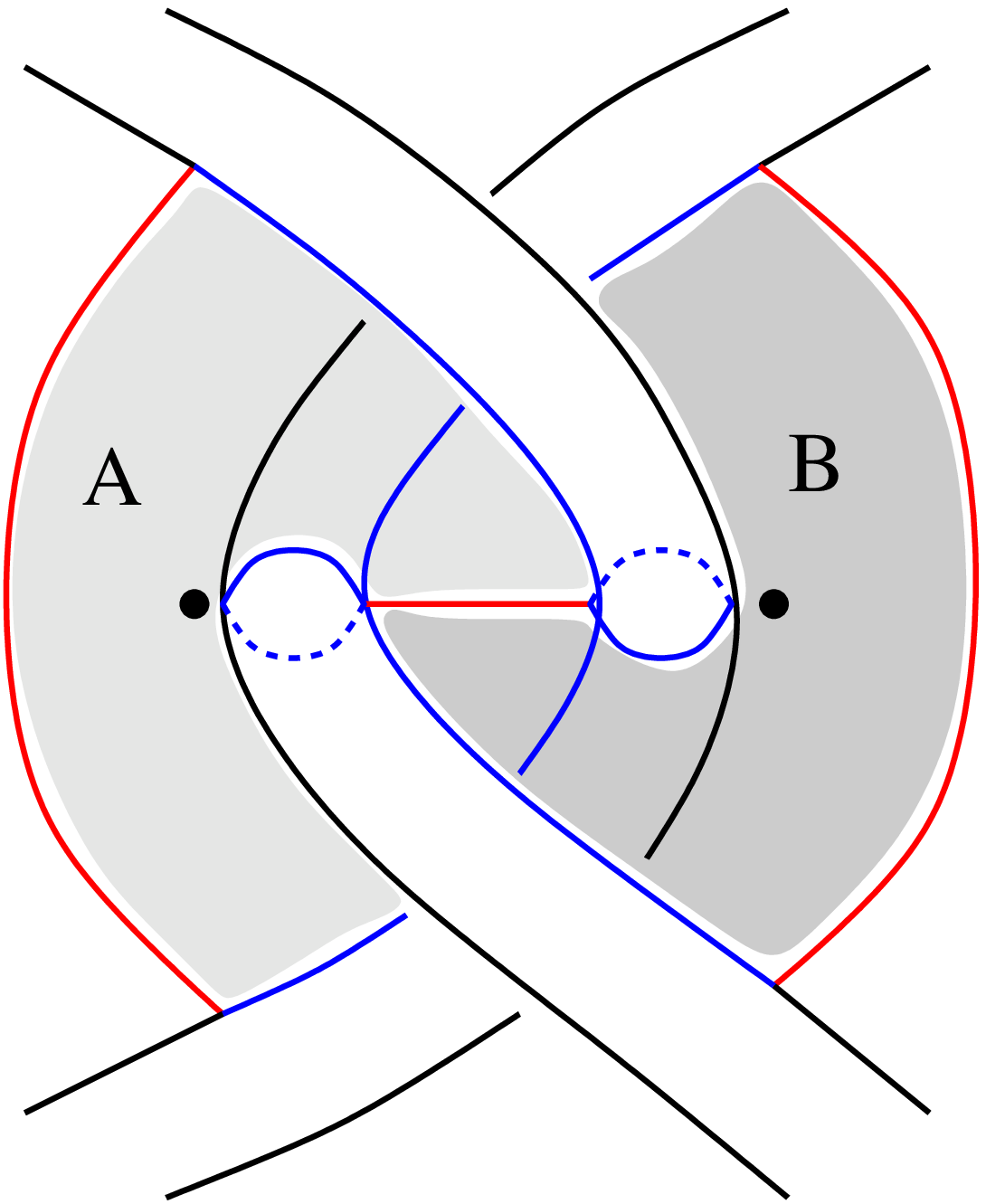}
\quad \quad
\includegraphics[width=80pt,clip]{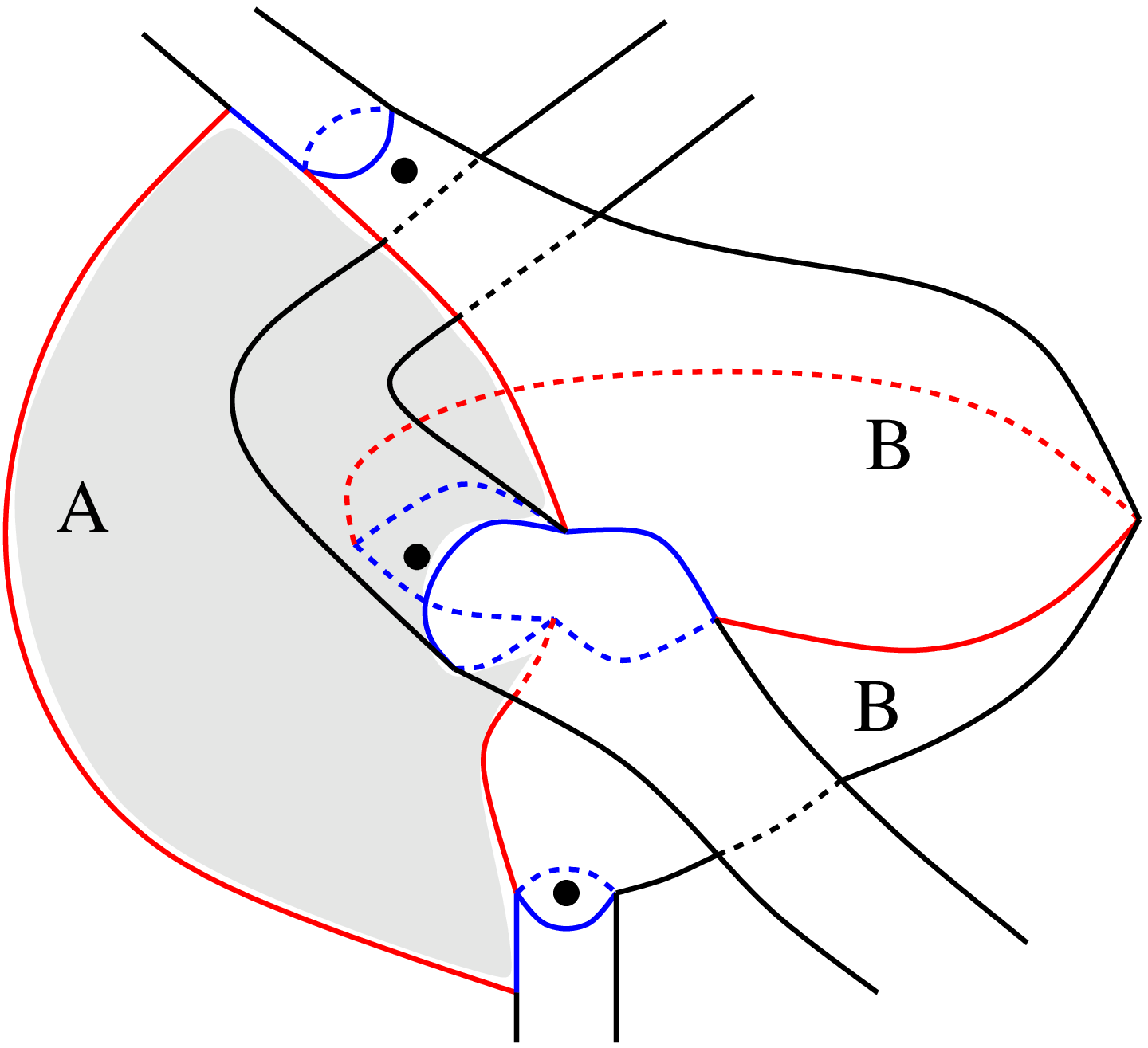}
\quad \quad
\includegraphics[width=60pt,clip]{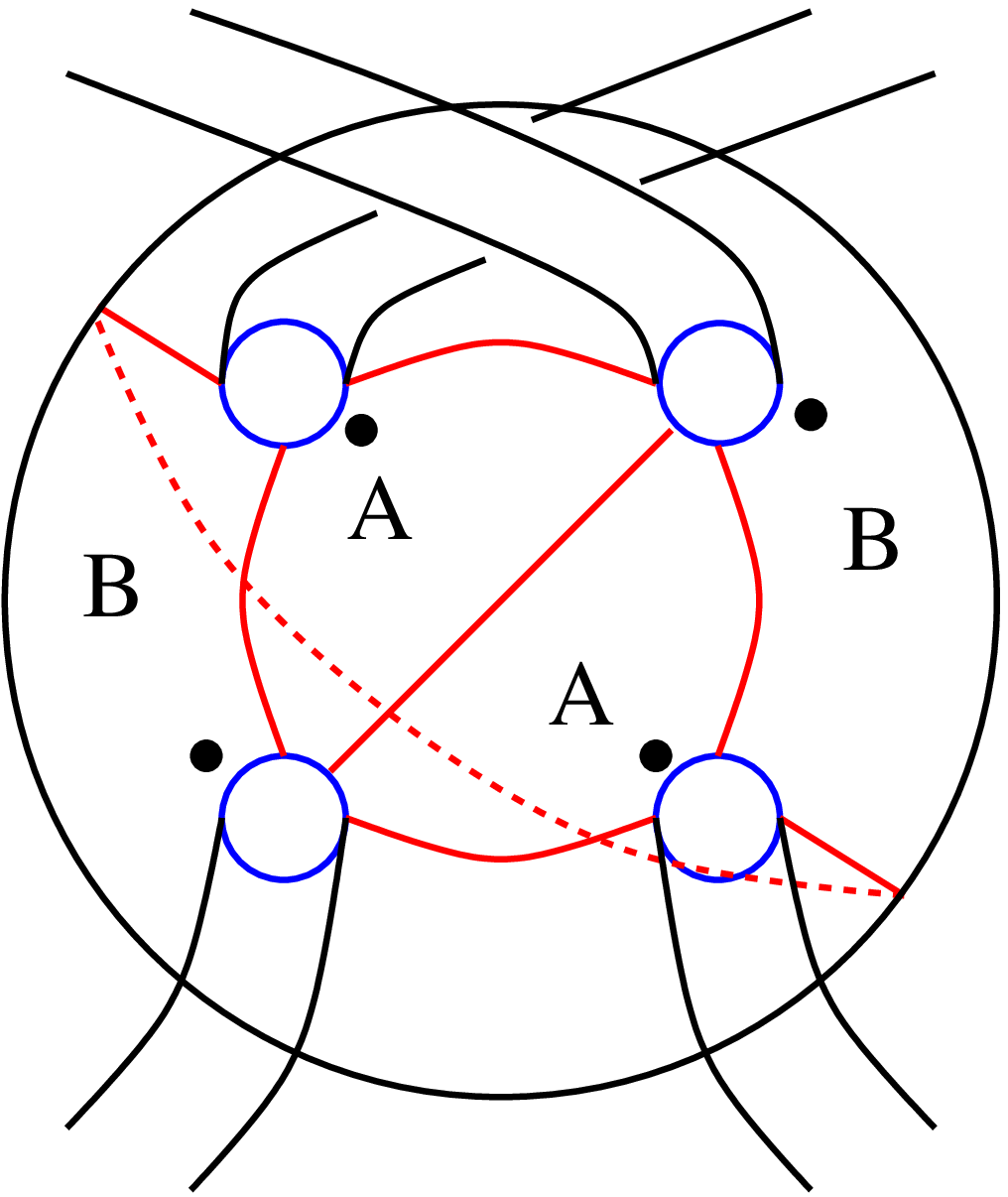}
\quad \quad
\includegraphics[width=60pt,clip]{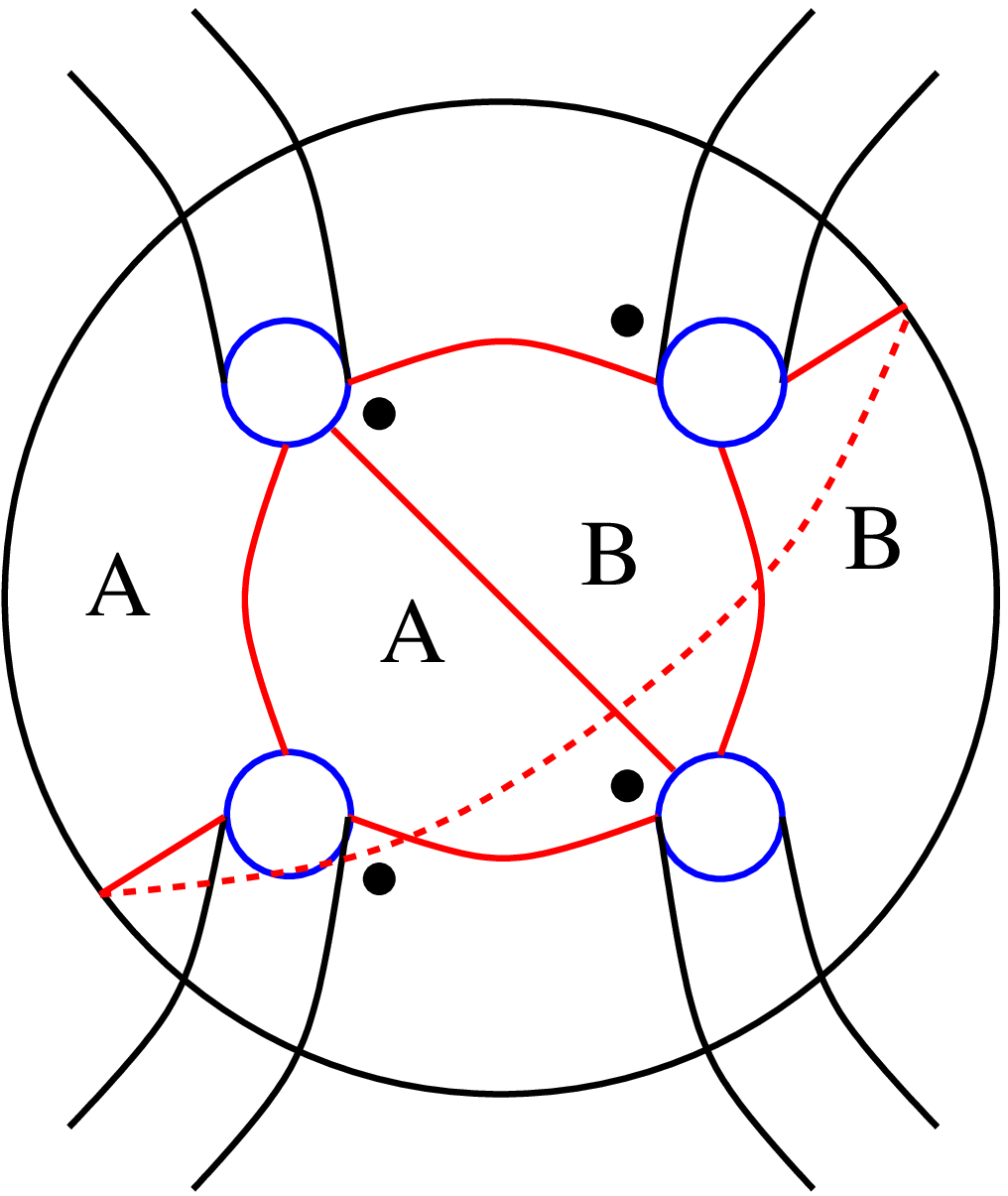}
\caption{Slice the full twist part.}
\label{fig:fulltwist}
\end{center}
\end{figure}

\begin{figure}[ht]
\begin{center}
\includegraphics[width=170pt, clip]{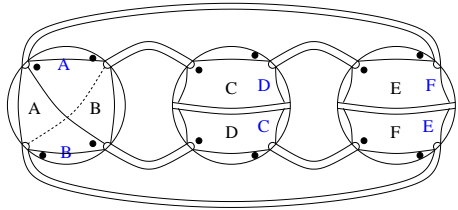}
\caption{The result of slicing the full twist and along the disks bounded by the trivial components.}
\label{fig:slice}
\end{center}
\end{figure}

\begin{figure}[ht]
\begin{center}
\includegraphics[width=150pt,clip]{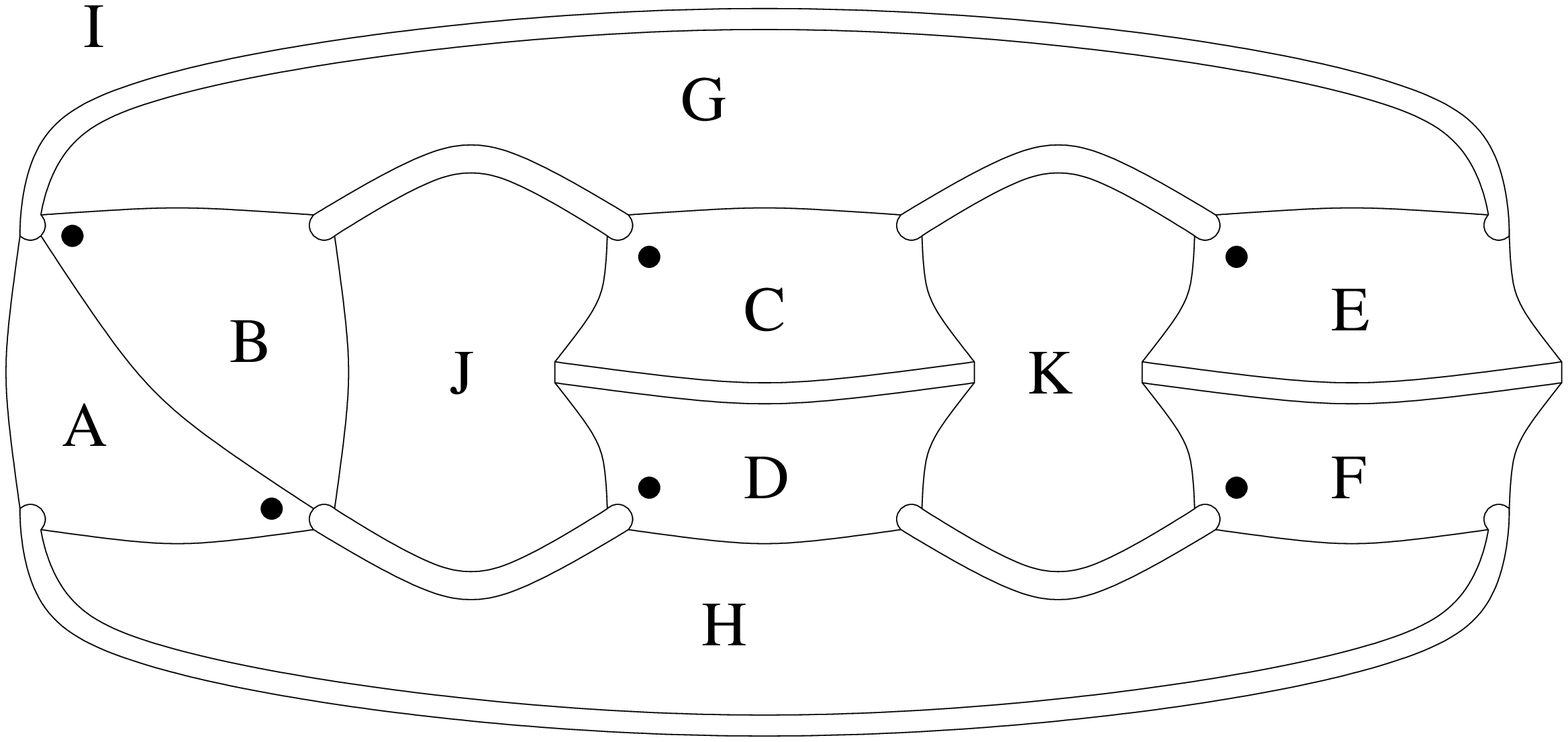}
\quad
\includegraphics[width=150pt,clip]{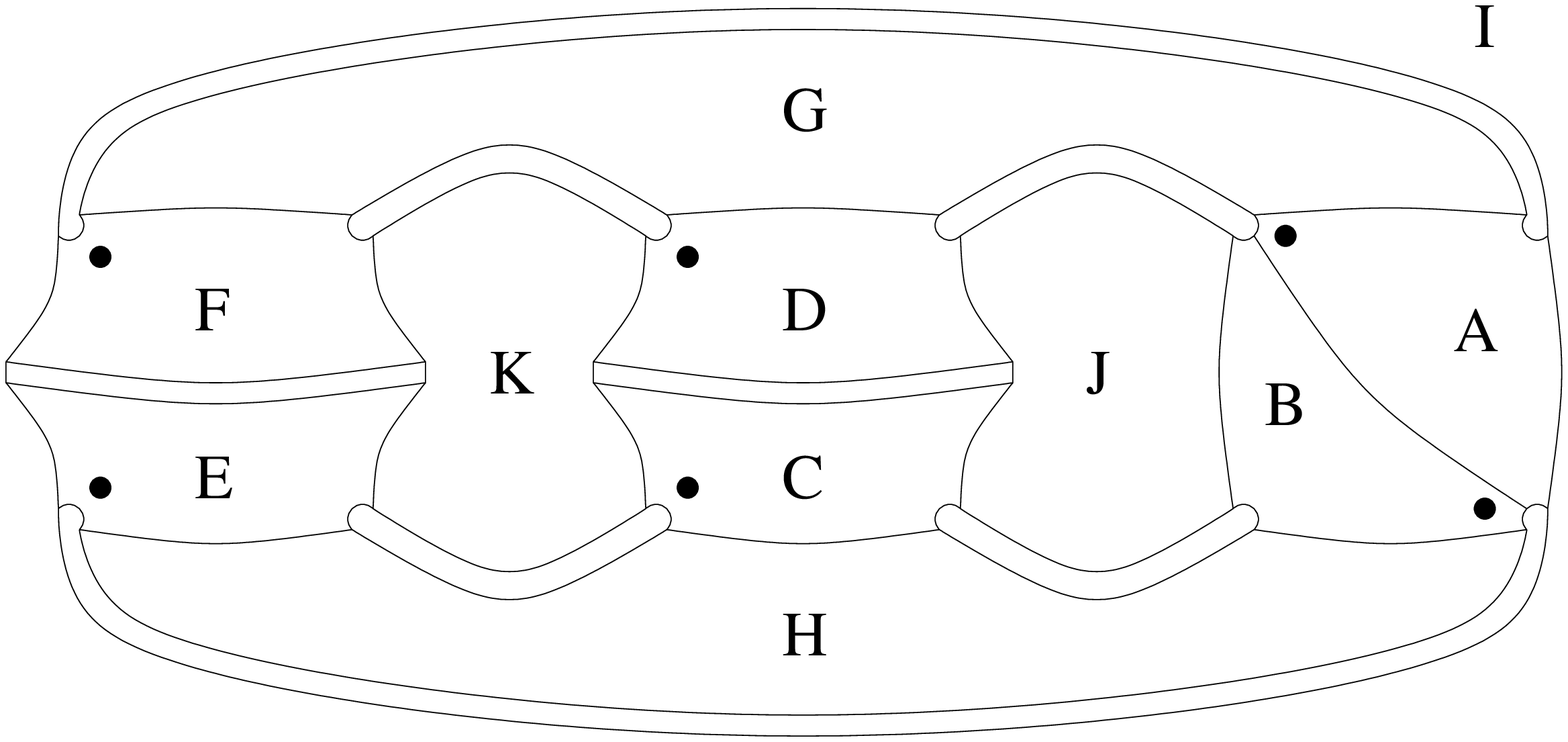}
\caption{The left is the `upper' side of Figure \ref{fig:slice} and 
the right is the `lower' side of Figure \ref{fig:slice}.
We are viewing these pictures from inside the 3-balls.}
\label{fig:upper-lower}
\end{center}
\end{figure}

\begin{figure}[hb]
\begin{center}
\includegraphics[width=120pt,clip]{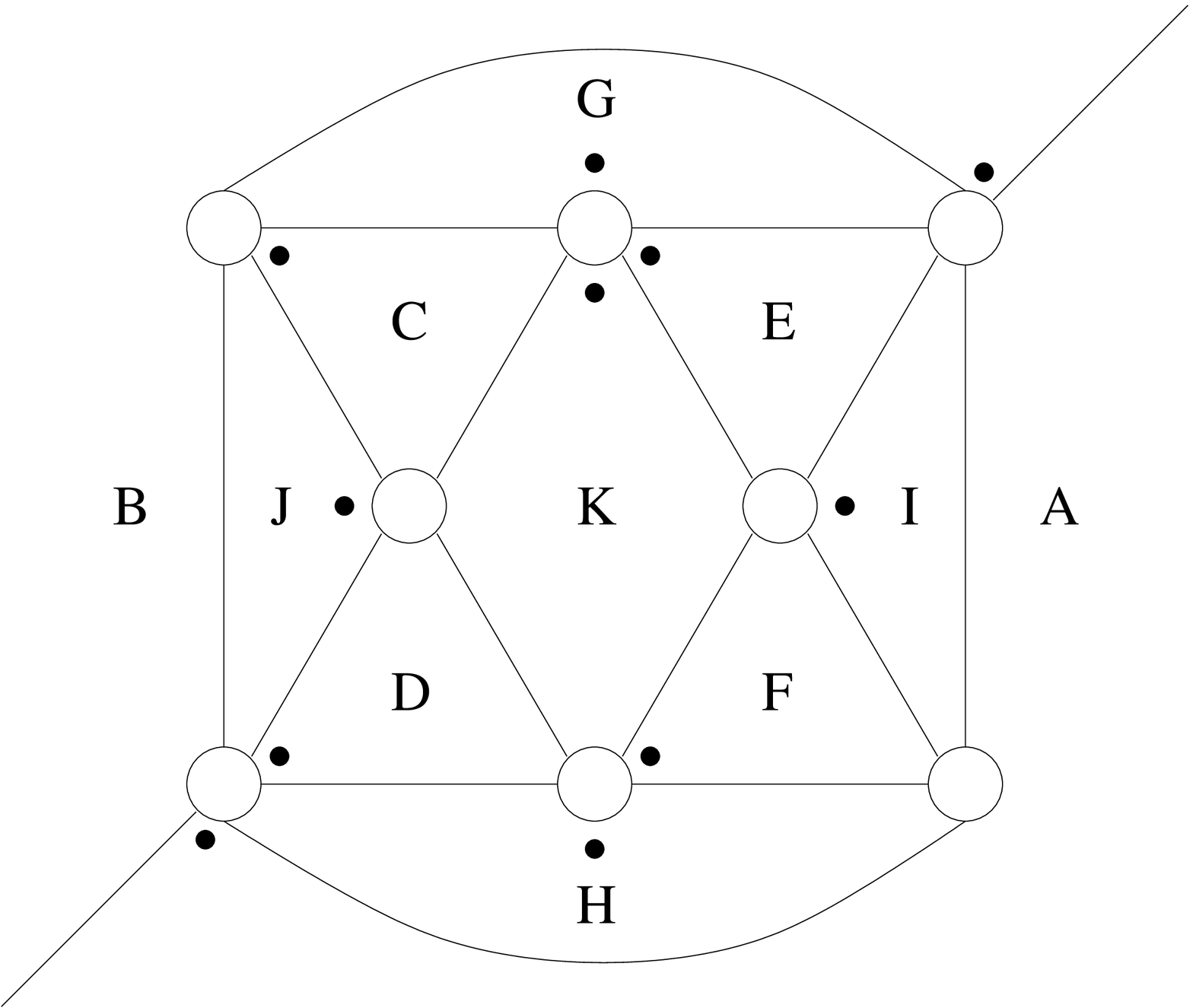}
\includegraphics[width=120pt,clip]{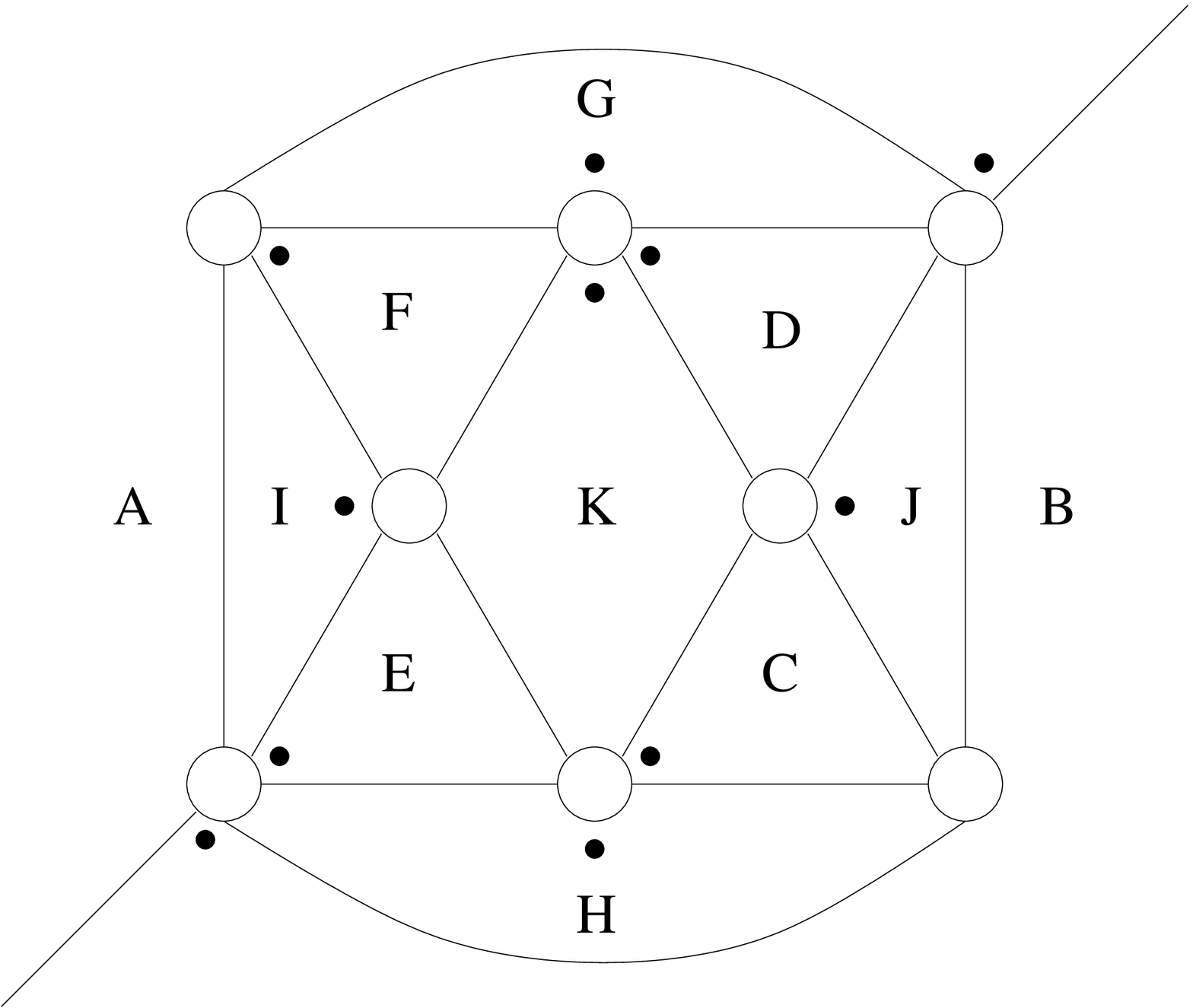}
\caption{Shrink `fragments of cusps' to small disks.}
\label{fig:upper-lower_2}
\end{center}
\end{figure}

Take a disk bounding each trivial link component and
cut $S^3\setminus L$ along these two disks.
We also slice the link complement at the full-twist part as shown in Figure \ref{fig:fulltwist}.
The result is Figure \ref{fig:slice}.
In that figure, faces are attached to each other so that the black dots in the faces coincide.
We then slice the manifold into two balls along the horizontal plane (Figure \ref{fig:upper-lower}).

Now, shrink fragments of boundary components to small disks.
This gives a (topological) ideal polyhedral decomposition of $S^3\setminus L$ into two ideal polyhedra
(Figure \ref{fig:upper-lower_2}).
In Figure \ref{fig:upper-lower_2}, we change the diagonal edges of the square made of 
faces $A$ and $I$. 
We also change the diagonal edge of the square formed by $B$ and $J$ leaving us with
two balls with graphs as shown in  Figure \ref{fig:upper-lower_3}. 
Glue these two 3-balls along the $K$-faces to arrive at 
the manifold of Figure \ref{fig:idealtriangulation}.
We ideally triangulate Figure \ref{fig:idealtriangulation} into ten ideal tetrahedra
as shown in Figure \ref{fig:hexa}.


\begin{figure}[ht]
\begin{center}
\includegraphics[width=100pt,clip]{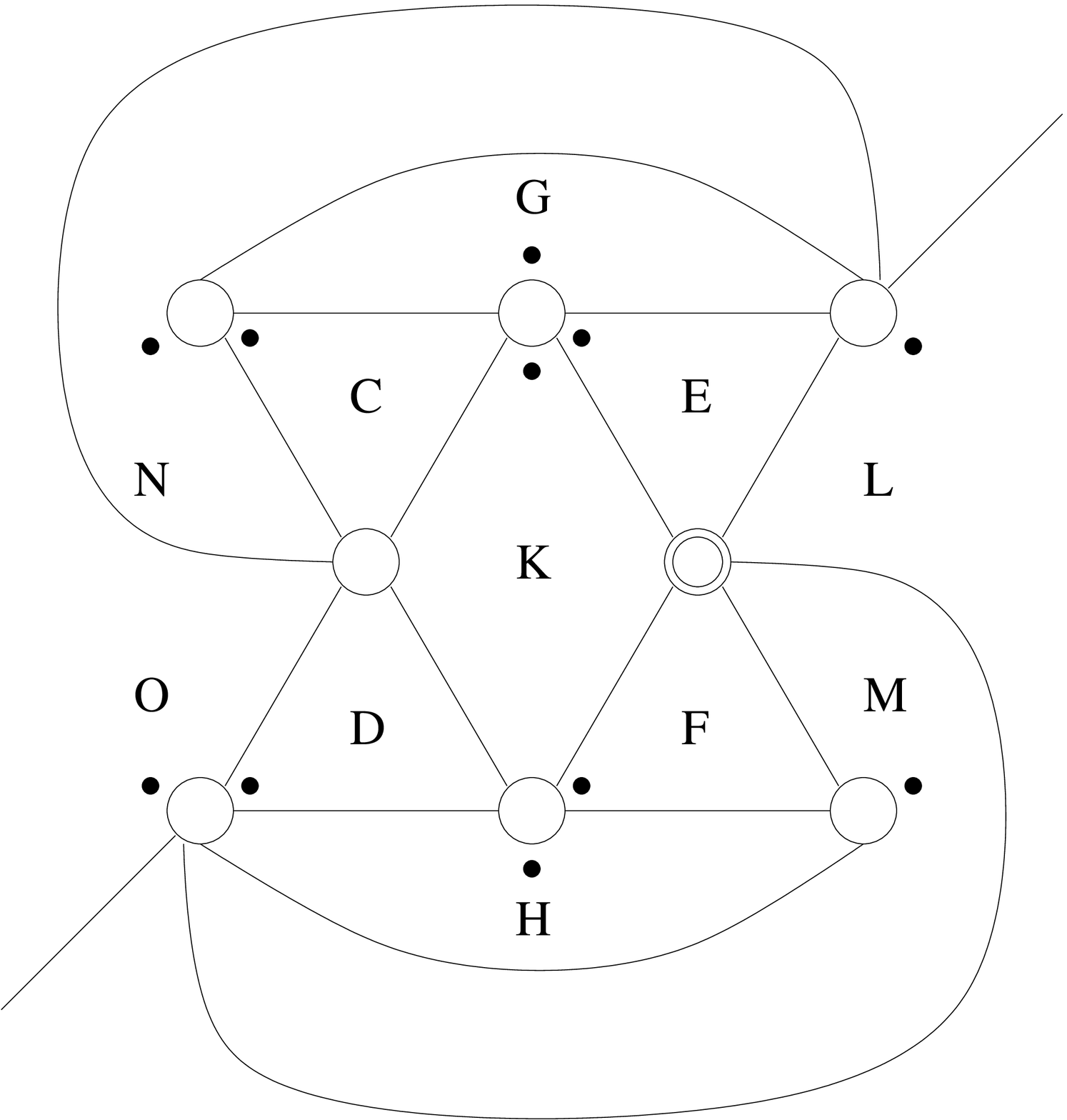}
\includegraphics[width=100pt,clip]{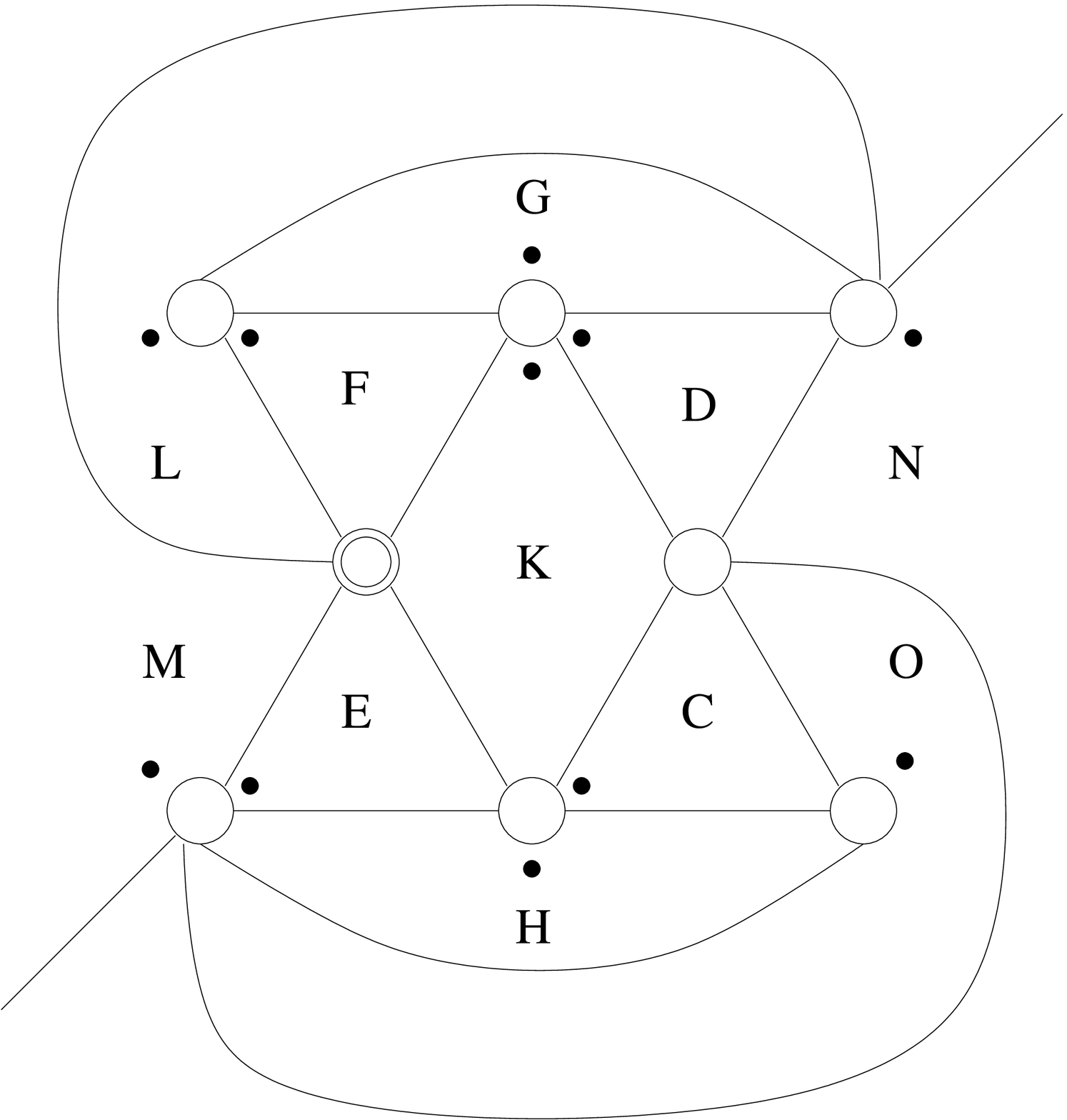}
\caption{Exchange the diagonals between $A$ and $I$ and also between $B$ and $J$.}
\label{fig:upper-lower_3}
\end{center}
\end{figure}

\begin{figure}[ht]
\begin{center}
\includegraphics[width=120pt,clip]{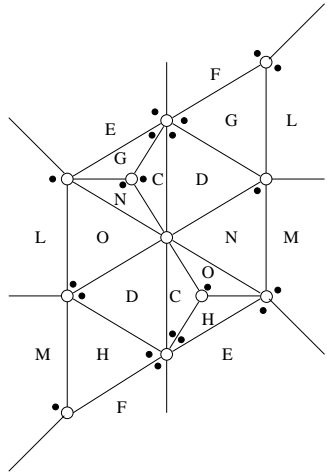}
\caption{An ideal triangulation of the link complement $S^3 \setminus L$.}
\label{fig:idealtriangulation}
\end{center}
\end{figure}


For an alternate point of view, we also note this decomposition into ten ideal tetrahedra is a refinement of the ideal polyhedral decomposition that Futer and Gu\'eritaud describe in \cite[Section 4]{fg:arborescent}.

We can easily observe that each 1-simplex of 
the ideal triangulation is the edge of $6$ ideal tetrahedra.
Therefore, if we give a regular ideal tetrahedral structure to each tetrahedron,
these ideal tetrahedra will satisfy the gluing equation around each 1-simplex.
We can see that the torus boundaries have the Euclidean structures shown in
Figures \ref{fig:originalcusp} and \ref{fig:cusp_trivial}.
Thus we obtain a complete hyperbolic structure on $S^3\setminus L$.

\begin{figure}[h]
\begin{center}
\hspace{30pt}
\includegraphics[width=80pt,clip]{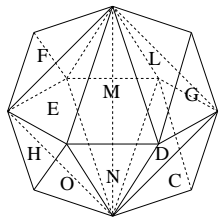}
\caption{The ideal polyhedron decomposed into $10$ ideal tetrahedra.}
\label{fig:hexa}
\end{center}
\end{figure}

\begin{figure}[ht]
\begin{center}
\includegraphics[width=170pt,clip]{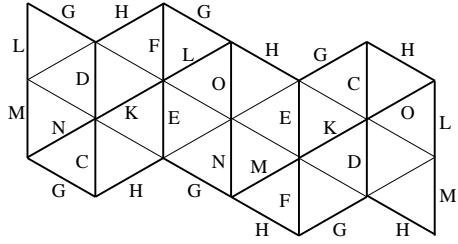}
\caption{Euclidean structure of torus boundary ($(-2,1,1)$ component).}
\label{fig:originalcusp}
\end{center}
\end{figure}

\begin{figure}[ht]
\begin{center}
\hspace{30pt}
\includegraphics[width=60pt,clip]{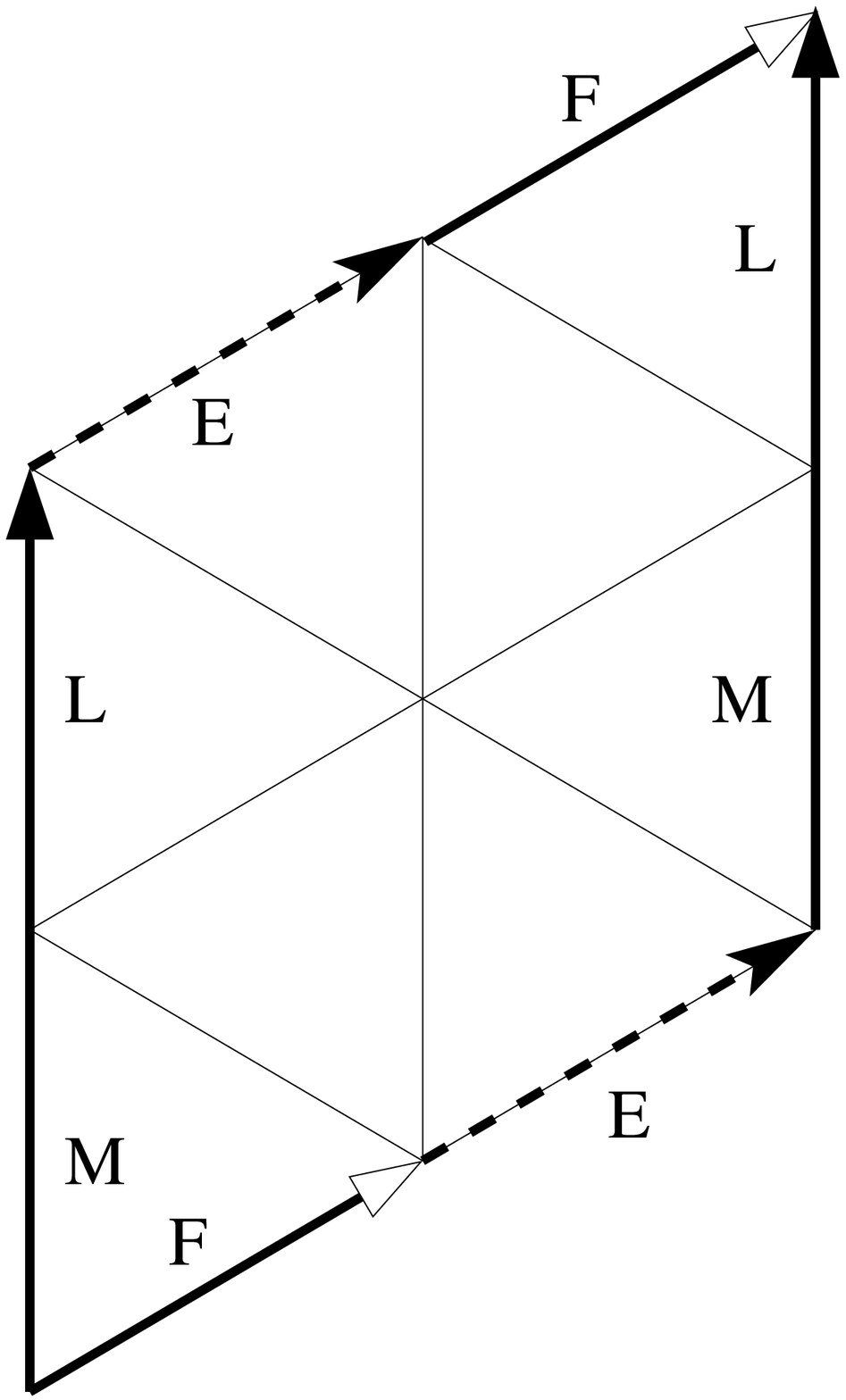}
\caption{Euclidean structure of torus boundary (trivial component).}
\label{fig:cusp_trivial}
\end{center}
\end{figure}

If we take a uniform cusp cross-section for each ideal tetrahedron as shown in Figure \ref{fig:regularideal}, the 
cusp cross-section has Euclidean structure as shown in Figure \ref{fig:originalcusp} and 
Figure \ref{fig:cusp_trivial} where the side of each triangle has Euclidean length 1.

\begin{figure}[ht]
\centerline{
\begin{picture}(0,0)%
\includegraphics{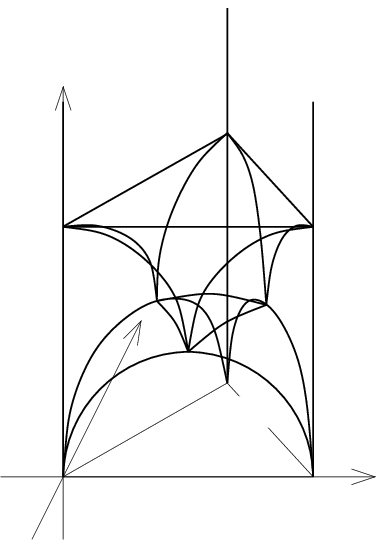}%
\end{picture}%
\setlength{\unitlength}{987sp}%
\begingroup\makeatletter\ifx\SetFigFont\undefined%
\gdef\SetFigFont#1#2#3#4#5{%
  \reset@font\fontsize{#1}{#2pt}%
  \fontfamily{#3}\fontseries{#4}\fontshape{#5}%
  \selectfont}%
\fi\endgroup%
\begin{picture}(7224,10245)(-1211,-6373)
\put(2401,-4261){\makebox(0,0)[lb]{\smash{{\SetFigFont{8}{9.6}{\familydefault}{\mddefault}{\updefault}$\frac{1+\sqrt{-3}}{2}$}}}}%
\put(226,-5611){\makebox(0,0)[lb]{\smash{{\SetFigFont{8}{9.6}{\familydefault}{\mddefault}{\updefault}$0$}}}}%
\put(4651,-5611){\makebox(0,0)[lb]{\smash{{\SetFigFont{8}{9.6}{\familydefault}{\mddefault}{\updefault}$1$}}}}%
\put(-449,-511){\makebox(0,0)[lb]{\smash{{\SetFigFont{8}{9.6}{\familydefault}{\mddefault}{\updefault}$1$}}}}%
\end{picture}%
}
\caption{A uniform cusp cross section.}
\label{fig:regularideal}
\end{figure}

We shall apply Theorem \ref{6theorem} to the above cusp cross-section.
The universal cover of a cusp is the Euclidean plane, and a basepoint on the torus lifts to  a lattice generated by two complex numbers. Every slope on the torus lifts to a primitive lattice point, where the length of the slope is equal to the Euclidean distance of the lattice point from the origin.
For the trivial component of $L$, Figure \ref{fig:cusp_trivial} shows that the lattice is generated by the complex numbers $2$ (corresponding to the longitude) and $\sqrt{3}i$ (corresponding to the meridian).
Thus the Euclidean length of slope $-1/k$ is equal to $|\sqrt{-3}-2k|=\sqrt{3+4k^2}$. 
For $k \geq 3$, (hence $p,q \geq 7$), the slope $-1/k$ is longer than 6.

At the cusp of the original knot, Figure \ref{fig:regularideal} shows that the lattice giving the Euclidean structure is generated by the complex numbers $2$ (corresponding to the meridian) and $-1+3\sqrt{3}i$ (corresponding to the slope $4/1$). After $-1/k$-surgery and $-1/l$-surgery on the trivial link components,
the slope $4/1$ becomes $4(k+l)+4=2(p+q)$ 
because the linking number with each of the trivial link components is equal to $2$.
The length of the slope $(2(p+q)m+n)/m$ is equal to $|2n+(-1+3\sqrt{3}i)m|=\sqrt{(2n-m)^2+27m^2}$.
Therefore, we can realize every surgery on $K$ via a Dehn filling of $L$ along slopes longer than $6$, except $(2(p+q)m+n)/m$ surgery when $(m,n)=(0,1), (1,-1), (1,0), (1,1), (1,2)$. 
By Theorem \ref{6theorem}, this means there are at most $5$ exceptional surgeries on $K$.
\end{proof}

\begin{proof}[Proof of Proposition \ref{propk2}]
In this case we have to replace the cusp cross-sections in order to enlarge one cusp 
so that the slope $-1/2$ on one of the trivial link components has length greater than $6$.
Then the other cusp cross-sections will become smaller.
We expand the cusp cross-section corresponding to $5$ half twists, by a factor of $\sqrt{2}$.
The length of slope $-1/2$ on that cusp is now equal to $\sqrt{2(3+4 \cdot 2^2)}=\sqrt{38} > 6$.
Then the other cusps are contracted by $1/\sqrt{2}$ (Figure \ref{fig:expandedcusp}).

\begin{figure}[ht]
\centerline{
\begin{picture}(0,0)%
\includegraphics{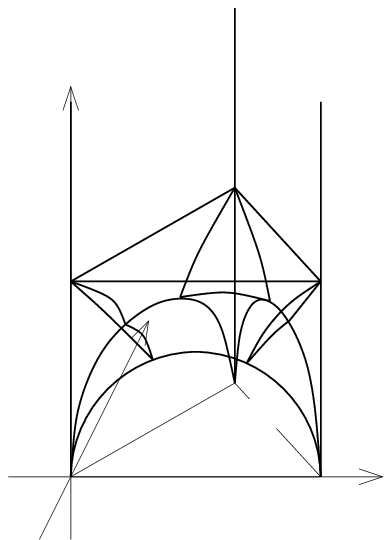}%
\end{picture}%
\setlength{\unitlength}{987sp}%
\begingroup\makeatletter\ifx\SetFigFont\undefined%
\gdef\SetFigFont#1#2#3#4#5{%
  \reset@font\fontsize{#1}{#2pt}%
  \fontfamily{#3}\fontseries{#4}\fontshape{#5}%
  \selectfont}%
\fi\endgroup%
\begin{picture}(7812,10245)(-1799,-6373)
\put(2401,-4261){\makebox(0,0)[lb]{\smash{{\SetFigFont{8}{9.6}{\familydefault}{\mddefault}{\updefault}$\frac{1+\sqrt{-3}}{2}$}}}}%
\put(226,-5611){\makebox(0,0)[lb]{\smash{{\SetFigFont{8}{9.6}{\rmdefault}{\mddefault}{\updefault}$0$}}}}%
\put(4651,-5611){\makebox(0,0)[lb]{\smash{{\SetFigFont{8}{9.6}{\familydefault}{\mddefault}{\updefault}$1$}}}}%
\put(-1799,-1561){\makebox(0,0)[lb]{\smash{{\SetFigFont{8}{9.6}{\familydefault}{\mddefault}{\updefault}$1/\sqrt{2}$}}}}%
\end{picture}%
}
\caption{Expanding one cusp cross-section at an ideal vertex.}
\label{fig:expandedcusp}
\end{figure}

After this modification of cusp neighborhoods, the Euclidean lattice for the other trivial link components is generated by the complex numbers $2/\sqrt{2}=\sqrt{2}$ (corresponding to the longitude) and 
$\sqrt{3}i/\sqrt{2}=\sqrt{6}i/2$  (corresponding to the meridian).
The Euclidean length of slope $-1/l$ is equal to $|\sqrt{6}i/2-\sqrt{2}l|=\sqrt{1.5+2l^2}$. 
So, for $l \geq 5$, the $-1/l$ slope is longer than $6$.

At the cusp of the original knot, the Euclidean lattice is generated by $\sqrt{2}$ (corresponding to the meridian)
and $-\sqrt{2}/2+3\sqrt{6}i/2$ (corresponding to the slope $4/1$).
After $-1/l$-surgery on the trivial link component,
the $4/1$-slope becomes $4(2+l)+4=2(5+q)$.
The length of the slope $(2(5+q)m+n)/m$ is equal to 
$|\sqrt{2}n+(-\sqrt{2}/2+3\sqrt{6}i/2)m|=\sqrt{(2n-m)^2/2+27m^2/2}$.
Thus, except for $(m,n)=(0,1), (1,-2), (1,-1), (1,0), (1,1), (1,2), (1,3)$, every surgery on $K$ can be realized by filling $L$ along slopes longer than $6$.
Therefore, there are at most $7$ exceptional surgeries.
\end{proof}

\subsection{\label{secpfive}%
Finite surgeries on the $(-2,5,q)$ pretzel}

In this subsection, we will prove
\begin{thm} \label{thmp5}
Let $q$ be odd with $q \geq 11$. The $(-2,5,q)$ pretzel knot admits no non-trivial finite surgery.
\end{thm}

\begin{proof}
By Proposition~\ref{propk2}, their are seven candidates for finite 
surgery. Lemma~\ref{lemfinbdy} eliminates $2q+10$ and $2q+12$ and Lemmas~\ref{lemeven} and \ref{lem2pqm1} rule out slopes $2q + 8$ and $2q+9$.

It remains to examine the slopes $2q+11$ and $2q+13$.
Since the boundary slopes of the $(-2,5,q)$ pretzel are
$0, 14, 15, \frac{q^2-q-5}{(q-3)/2}$, $2q+10$, and $2q+12$,
(see~\cite{ho}) the norm is 
\begin{eqnarray*}
\| \gamma \| & =&  2 {[} a_1 \Delta(\gamma, 0) +
a_{2} \Delta(\gamma, 14) + a_{3} \Delta(\gamma, 15) + 
a_{4} \Delta(\gamma, \frac{q^2-q-5}{(q-3)/2}) \\
&& \mbox{ } +
a_{5} \Delta(\gamma, 2q+10) + a_{6} \Delta(\gamma, 2q+12) {]} 
\end{eqnarray*}
where the $a_i$ are non-negative integers.

Thus, the norms of $\frac10$, $2q+10$, and $2q+11$ are
\begin{eqnarray*}
\|\frac10 \| & =&   S = 2 {[} a_1 + a_2 + a_3 + \frac{q-3}{2} a_4  + a_5 +
a_6 {]}, \\
\| 2q+10 \| & =&  2 {[} a_1 (2q+10) +
a_{2} (2q - 4) + a_{3} (2q-5) + 
a_{4} (3q-10) + 2 a_{6} {]},  \mbox{ and } \\
\| 2q+11 \| & =& 2 {[} a_1 (2q+11) +
a_{2} (2q - 3) + a_{3} (2q-4) + 
a_{4} (7q-23)/2 + a_5 + a_{6} {]}. 
\end{eqnarray*}

But then, 
$$\| 2q+11 \| - \|2q+10 \| = 2{[} a_1 + a_2 + a_3 + \frac{q-3}{2} a_4  + a_5 - a_6 {]}. $$
In other words, $\| 2q+11 \| = \| 2q+10 \| + S - 4 a_6$. 

Suppose $a_i = 0$ for $i \leq 4$. 
Then $S = 2(a_5 + a_6)$ and $\|2q+11 \| = 2(a_5 + a_6) = S$. 
As in~\cite{d1}, this would imply that there is a non-integral 
boundary slope $r$ with $|2q+11-r|<1$. As there is no such
$r$, we conclude that 
$\exists i \leq 4$ with $a_i>0$. Then
$\|2q+10\| \geq 2(2q-5 + 2a_6)$ and
$\|2q+11\| \geq S + 2(2q-5) > S+8$ so that $2q+11$ surgery is not
finite. 

A similar argument shows that $\|2q+13\|-\|2q+10\|=3S-4a_6$.
Hence we have $\|2q+13\|\geq 3S+2(2q-5)>S+8$ so that $2q+13$ is not finite.
\end{proof}

\subsection{\label{secpseven}%
Finite surgeries on the $(-2,p,q)$ pretzel ($7 \leq p \leq q$)}

In this subsection, we will prove
\begin{thm} \label{thmp7}
Let $p$ and $q$ be odd with $7 \leq p \leq q$. The $(-2,p,q)$ pretzel knot admits no non-trivial finite surgery.
\end{thm}

\begin{proof}
By Proposition~\ref{propk1}, their are five candidates for finite 
surgery. Lemma~\ref{lemfinbdy} eliminates $2(p+q)$ and $2(p+q)+2$ and Lemma~\ref{lem2pqm1} rules out the slope $2(p+q)-1$.

So, the only candidate for a non-trivial finite surgery is $2(p+q)+1$
and, by Lemma~\ref{lem2pq1}, we can assume either
$p = q= 9$ or else $p=7$ and $7 \leq q \leq 19$. However, an argument
similar to that used in the previous subsection for 
$2q+11$ surgery on the $(-2,5,q)$ pretzel knot shows that these remaining
eight cases also do not lead to a finite surgery.
\end{proof}

\end{document}